\documentclass[reqno]{amsart}  
\theoremstyle{plain}
\newtheorem{theorem}{Theorem}[section]
\newtheorem{lemma}[theorem]{Lemma}
\newtheorem{remark}[theorem]{Remark}
\newtheorem{proposition}[theorem]{Proposition}

\newtheorem{corollary}[theorem]{Corollary}

\numberwithin{equation}{section}
\allowdisplaybreaks[1]

\theoremstyle{definition}

\theoremstyle{remark}

\newcommand{\bA}{{\mathbf A}}
\newcommand{\bU}{{\mathbf U}}
\newcommand{\bZ}{{\mathbf Z}}

\newcommand{\blam}{{\boldsymbol \lambda}}
\newcommand{\lam}{\lambda}
\newcommand{\bzeta}{{\boldsymbol \zeta}}

\newcommand{\cC}{{\mathcal C}}

\newcommand{\cF}{{\mathcal F}}
\newcommand{\cH}{{\mathcal H}}
\newcommand{\cL}{{\mathcal L}}
\newcommand{\cM}{{\mathcal M}}

\newcommand{\cS}{{\mathcal S}}

\newcommand{\cU}{{\mathcal U}}
\newcommand{\cX}{{\mathcal X}}
\newcommand{\cY}{{\mathcal Y}}
\newcommand{\cO}{{\mathcal O}}

\newcommand{\DB}{{\mathcal H}^\tau(K_S)}

\newcommand{\sbm}[1]{\left[\begin{smallmatrix} #1
		\end{smallmatrix}\right]}

\begin{document}

\title[Multipliers on the Fock space]{Schur-class multipliers on the
Fock space: de Branges-Rovnyak reproducing kernel spaces and
transfer-function realizations}
\author[J. A. Ball]{Joseph A. Ball}
\address{Department of Mathematics,
Virginia Tech,
Blacksburg, VA 24061-0123, USA}
\email{ball@math.vt.edu}
\author[V. Bolotnikov]{Vladimir Bolotnikov}
\address{Department of Mathematics,
The College of William and Mary,
Williamsburg VA 23187-8795, USA}
\email{vladi@math.wm.edu}
\author[Q. Fang]{Quanlei Fang}
\address{Department of Mathematics,
Virginia Tech,
Blacksburg, VA 24061-0123, USA}
\email{qlfang@math.vt.edu}

\begin{abstract}
We introduce and study a Fock-space noncommutative analogue of reproducing
kernel Hilbert spaces of de Branges-Rovnyak type.
Results include: use of the de Branges-Rovnyak space $\cH(K_{S})$ as
the state space for the unique (up to unitary equivalence)
observable, coisometric transfer-function realization of the
Schur-class multiplier $S$, realization-theoretic characterization of
inner Schur-class multipliers, and a calculus for obtaining a
realization for an inner multiplier with prescribed left zero-structure.
In contrast with the parallel theory for the Arveson space on the
unit ball ${\mathbb B}^{d} \subset {\mathbb C}^{d}$ (which can be
viewed as the symmetrized version of the Fock space used here), the results
here are much more in line with the classical univariate case, with
the extra ingredient of the existence of all results having both a
``left'' and a ``right'' version.

\end{abstract}


\maketitle

\centerline{\em Dedicated to the memory of Tiberiu Constantinescu}

\section{Introduction}  \label{S:Intro}
\setcounter{equation}{0}

Recently there has been much interest and an evolving theory of
    noncommutative function theory and associated multivariable
    operator theory and multidimensional system theory with evolution
    along a free semigroup; we mention \cite{A-KV, KV-V, BB-noncomint, BGM1,
    BGM2, CJ, HMcCV, MS-Annalen, MS-Schur, PopescuNF1, PopescuNF2,
    Popescu-Nehari, Popescu-Memoir}.
    A central player in many of these developments is the
    noncommutative Schur class consisting of formal power series in a
    set of noncommuting indeterminates which define contractive
    multipliers between (unsymmetrized) vector-valued Fock spaces; such
    Schur-class functions play the role of the characteristic function
    for the Popescu analogue for a row contraction of the
    Sz.-Nagy-Foia\c s model theory for a single contraction operator
    (see \cite{PopescuNF2, Cuntz-scat}). For the classical (univariate) case,
    there is an approach to operator-model theory complementary to the
    Sz.-Nagy-Foia\c s approach which emphasizes constructions with
reproducing kernel
    Hilbert spaces over the unit disk rather than the geometry of the
    unitary dilation space of a contraction operator.  Our purpose here
    is to flesh out the ingredients of this approach for the Fock space
    setting.  The appropriate noncommutative multivariable version of a
    reproducing kernel Hilbert space has already been worked out in
    \cite{NFRKHS} and certain other relevant background material
    appears in \cite{BBF1}.  Unlike the work in some of the
    papers mentioned above, specifically
    \cite{A-KV, ariaspopescu,  BB-noncomint, BGM2, CJ,
    davidsonpitts, HMcCV, KV-V, MS-Annalen, MS-Schur, Popescu-Nehari},
    we shall deal with formal power series with operator coefficients
    as parts of some formal structure
    (e.g., as inducing multiplication operators between two Hilbert
    spaces whose elements are formal power series with vector
    coefficients) rather than as themselves functions on some
    collection of noncommutative operator-tuples.
    Before discussing the precise
    noncommutative results which we present here, we review the
    corresponding classical versions of the results.

For $\cU$ and $\cY$ two Hilbert spaces, let $\cL(\cU, \cY)$ denote
the space of bounded linear operators between $\cU$ and $\cY$.
We also let $H^2_{\cU}({\mathbb D})$ be the standard Hardy
space of the $\cU$-valued holomorphic functions on the unit disk ${\mathbb
D}$. By the classical Schur class ${\mathcal S}(\cU, \cY)$ we mean the
set of $\cL(\cU, \cY)$-valued functions holomorphic on the unit
disk ${\mathbb D}$ with values $S(\lam)$ having norm at most $1$ for
each $\lam \in {\mathbb D}$.  There are several equivalent
characterizations of the class ${\mathcal S}(\cU, \cY)$; for
convenience, we list some in the following theorem.

\begin{theorem} \label{T:I1}
       Let $S$ be an $\cL(\cU, \cY)$-valued function defined on the unit
       disk ${\mathbb D}$.  Then the following are equivalent:
       \begin{enumerate}
	\item $S \in \cS(\cU, \cY)$, i.e., $S$ is analytic on ${\mathbb D}$
	with contractive values in $\cL(\cU, \cY)$.
\item The multiplication operator $M_{S} \colon
            f(z) \mapsto S(z) \cdot f(z)$ is a contraction from
            $H^{2}_{\cU}({\mathbb D})$ into $H^{2}_{\cY}({\mathbb D})$.
	\item The kernel
	$$ K_{S}(\lam, \zeta) := \frac{ I_{\cY} - S(\lam) S(\zeta)^{*}}{1 -
	\lam \overline{\zeta}}
	$$
	is positive on ${\mathbb D} \times {\mathbb D}$, i.e.,
	there exists an auxiliary Hilbert space $\cX$ and a function $H
	\colon {\mathbb D} \to \cL(\cX, \cY)$ such that
	\begin{equation} \label{KSfact}
	  K_{S}(\lam, \zeta) = H(\lam) H(\zeta)^{*} \quad\text{for all}\quad
          \lam, \zeta \in {\mathbb D}.
	\end{equation}
	\item There exists a Hilbert space $\cX$ and a unitary connection
	operator (or colligation) $\bU$ of the form
	\begin{equation} \label{I:colligation}
	  \bU = \begin{bmatrix} A & B \\ C & D \end{bmatrix} \colon
	  \begin{bmatrix} \cX  \\ \cU \end{bmatrix} \to \begin{bmatrix} \cX
	  \\ \cY \end{bmatrix}
	\end{equation}
	so that $S(\lam)$ can be realized in the form
	\begin{equation}  \label{I:realization}
	 S(\lam) = D + \lam C (I_{\cX} - \lam A)^{-1} B.
	\end{equation}

	\item There exists a Hilbert space $\cX$ and a contractive
	connecting operator $\bU$ of the form \eqref{I:colligation} so that
	\eqref{I:realization} holds.
	\end{enumerate}
	\end{theorem}
A pair $(C,A)$ is called {\em an output pair} if  $C \in \cL(\cX, \cY)$
and  $A \in \cL(\cX, \cX)$. An output pair $(C,A)$ is called {\em
contractive} if $A^{*} A + C^{*} C \le I_{\cX}$, {\em isometric} if
$A^{*} A + C^{*} C = I_{\cX}$ and {\em observable} if ${\displaystyle
\bigcap_{n=0}^{\infty} \operatorname{Ker}\, C A^{n} =\{0\}}$.
We shall say that the realization \eqref{I:realization} of
$S(\lambda)$ is {\em observable} if the output pair $(C,A)$ occurring in
\eqref{I:realization} is observable.
Furthermore,
with an output contractive pair $(C,A)$, one can associate the positive
kernel
\begin{equation} \label{I:defKCA}
K_{C,A}(\lam, \zeta) = C(I - \lam A)^{-1} (I - \overline{\zeta}
A^{*})^{-1} C^{*}
\end{equation}
which is (as it is readily seen) defined on ${\mathbb D}\times{\mathbb D}$.

          As also remarked in \cite{BBF2}, the coisometric version
	of (4) $\Longrightarrow$ (2) is particularly transparent,
	since in this case a simple computation shows that then
          \eqref{KSfact} holds with $H(\lam) = C (I - \lam A)^{-1}$, i.e.,
          $K_{S}(\lam, \zeta) =K_{C,A}(\lam, \zeta)$. We have the following
          sort of converse of these observations.

	\begin{theorem} \label{T:I2}
	    \begin{enumerate}
		\item
	    Suppose that $S\in \cS(\cU, \cY)$ and that $(C,
	    A)$ is an observable, contractive output-pair
	   of operators  such that
	    \begin{equation} \label{I:KS=KCA}
	      K_{S}(\lam, \zeta) = K_{C,A}(\lam, \zeta).
	    \end{equation}
	  Then there is a unique choice of $B \colon \cU \to \cX$ so that
	  $\bU = \sbm{ A & B \\ C & S(0) }$ is coisometric and $\bU$
	  provides a realization for $S$:
	  $S(\lam) = S(0) + \lam C (I - \lam A)^{-1}B$.

	  \item Suppose that we are given only an observable,
	  contractive output-pair of operators $(C, A)$ as above.
	  Then there is a choice of an input space $\cU$ and a Schur
	  multiplier $S \in {\mathcal S}(\cU, \cY)$ so that
	  \eqref{I:KS=KCA} holds.
	  \end{enumerate}
	\end{theorem}

	As we see from Theorem \ref{T:I1}, for any Schur-class function $S
	\in \cS(\cU, \cY)$, we can associate
	the positive kernel
	$K_{S}(\lam, \zeta)$ and therefore also by Aronszajn's
	construction the reproducing kernel
	Hilbert space $\cH(K_{S})$; this space is called the de
	Branges-Rovnyak space associated with $S$.  It turns out that 
any observable
	coisometric realization $\bU$ for $S$ is unitarily equivalent to a
	certain canonical functional-model realization.

	\begin{theorem} \label{T:I3}  Let $S\in \cS(\cU,\cY)$. Then the
          operator
          $$
          \bU_{\text{dBR}} = \begin{bmatrix} A_{\text{dBR}} &
          B_{\text{dBR}} \\ C_{\text{dBR}} & D_{\text{dBR}} \end{bmatrix}
          \colon \begin{bmatrix} \cH(K_{S}) \\ \cU \end{bmatrix} \to
\begin{bmatrix}
          \cH(K_{S}) \\ \cY \end{bmatrix}
          $$
with the entries given by
	\begin{align*}
	 A_{\text{dBR}} \colon f(\lam) \to \frac{f(\lam) - f(0)}{\lam},
	 & \qquad B_{\text{dBR}} \colon u \to
	 \frac{S(\lam) - S(0)}{\lam} u, \\
	  C_{\text{dBR}} \colon f \to f(0),
	 & \qquad D_{\text{dBR}} \colon u \to S(0) u
	\end{align*}
provides an  observable and coisometric realization
\begin{equation}
S(\lam) = D_{\text{dBR}}  + \lam C_{\text{dBR}}  (I_{\cH(K_{S})} - \lam
A_{\text{dBR}} )^{-1} B_{\text{dBR}}.
\label{debr}
\end{equation}
Moreover, any other observable coisometric realization of $S$ is
unitarily equivalent to \eqref{debr}.
\end{theorem}
Let us say that a Schur function $S \in \cS(\cU, \cY)$ is {\em inner}
if the associated multiplication operator $M_{S} \colon
H^{2}_{\cU}({\mathbb D}) \to H^{2}_{\cY}({\mathbb D})$ is a
partial isometry. Equivalently, $S \in \cS(\cU, \cY)$ and the almost
everywhere existing boundary value function $S(\zeta) = \lim_{r
\uparrow 1}S(r \zeta)$ is a partial isometry for almost all $\zeta \in
{\mathbb T}$.
The following characterization of inner
functions in terms  of realizations is well known (see
\cite{dbr1, dbr2}).

      \begin{theorem}\label{T:I5}  A Schur multiplier $S \in \cS(\cU,
      \cY)$ is inner if and only if its essentially unique observable,
      coisometric realization of the form
      \eqref{I:realization} is such that $A$ is strongly stable, i.e.,
      \begin{equation} \label{I:stable}
\lim_{n \to \infty} \| A^{n} x \| = 0 \text{ for all } x \in\cX.
      \end{equation}
      \end{theorem}

      Inner functions come up  in the representation of shift-invariant
      subspaces of $H^{2}_{\cY}$ as in the Beurling-Lax theorem.
      The following version of the Beurling-Lax theorem first identifies
      any shift-invariant subspace as the set of solutions of a
      collection of homogeneous interpolation conditions and then
      obtains a realization for the Beurling-Lax representer in terms of
      the data set for the homogeneous interpolation problem.
       The finite-dimensional version of this result can be found in
\cite[Chapter 14]{BGR} while the
      details of the general case appear in \cite{BallRaney}.
      We let $M_{\lam}$ denote the shift operator
      $$
      M_{\lam} \colon f(\lam) \to \lam f(\lam) \quad \text{for} \quad f \in
      H^{2}_{\cY}({\mathbb D})
      $$
and given a contractive pair $(C,A)$ we let
\begin{equation}  \label{defMCA}
{\mathcal M}_{A^{*},C^{*}} = \{ f \in H^{2}_{\cY}({\mathbb D})
\colon \;  (C^{*}f)^{\wedge L}(A^{*}) = 0 \}
\end{equation}
     where we have set
     $$
       (C^{*}f)^{\wedge L}(A^{*}): = \sum_{n=0}^{\infty} A^{*n} C^{*}
       f_{n} \quad\text{if} \quad f(\lam) = \sum_{n=0}^{\infty} f_{n}
\lam^{n}     \in H^{2}_{\cY}({\mathbb D}).
     $$
      \begin{theorem} \label{T:I6}  \begin{enumerate}
          \item Suppose that ${\mathcal M}$ is a subspace of
          $H^{2}_{\cY}({\mathbb D})$ which is $M_{\lam}$-invariant.
          Then there is an isometric pair $(C,A)$ such that $A$ is strongly
          stable (i.e., \eqref{I:stable} holds) and such that
          ${\mathcal M}={\mathcal M}_{A^{*},C^{*}}$.
     \item Suppose that the shift-invariant subspace ${\mathcal M}
     \subset H^{2}_{\cY}({\mathbb D})$ has the representation ${\mathcal
     M} = {\mathcal M}_{A^{*},C^{*}}$ as in \eqref{defMCA} where $(C,A)$
     is an isometric pair with $A$ strongly stable.  Choose an input
     space ${\mathcal U}$ and operators
     $B \colon \cU \to \cX$ and $D \colon \cU \to \cY$ so that
     $$
      \bU  = \begin{bmatrix} A & B \\ C & D \end{bmatrix} \colon
      \begin{bmatrix} \cX \\ \cU \end{bmatrix} \to \begin{bmatrix} \cX
      \\ \cY \end{bmatrix}
     $$
     is unitary.  Then the function
     $S(\lam) = D + \lam C (I_{\cX} - \lam A)^{-1} B$
     is inner (i.e., $M_{S}$ is isometric) and is a
     Beurling-Lax representer for ${\mathcal M}$:
     $$  S \cdot H^{2}_{\cU}({\mathbb D}) = {\mathcal M}_{A^{*},C^{*}}.
     $$
     \end{enumerate}
     \end{theorem}

     Our goal here is to obtain noncommutative analogues of these
     results, where the classical Schur class is replaced by the
     noncommutative Schur class of contractive multipliers between Fock
     spaces of formal power series in noncommuting indeterminates and
     where the classical reproducing kernel Hilbert spaces become the
     noncommutative formal reproducing kernel Hilbert spaces introduced
     in \cite{NFRKHS}.
     Let $z = (z_{1}, \dots, z_{d})$ and $w = (w_{1}, \dots, w_{d})$ be
     two sets of noncommuting indeterminates.  We let $\cF_{d}$ denote
     the free semigroup generated by the $d$ letters $\{1, \dots, d\}$.
     A generic element of $\cF_{d}$ is a word $w$ equal to a string of
     letters
     \begin{equation} \label{word}
      \alpha = i_{N} \cdots i_{1}\quad \text{where}\quad i_{k} \in \{1,
\dots, d\} \;    \text{ for } \; k=1, \dots, N.
      \end{equation}
      Given two words $\alpha$ and $\beta$ with $\alpha$ as in
      \eqref{word} and $\beta$ of
      the form $\beta = j_{N'} \cdots j_{1}$, say, the product $\alpha
      \beta$ is defined by concatenation:
      $$
        \alpha \beta = i_{N} \cdots i_{1} j_{N'} \cdots j_{1} \in \cF_{d}.
      $$
      The unit element of $\cF_{d}$ is the {\em empty word} denoted by
$\emptyset$.
      For $\alpha$ a word of the form \eqref{word}, we let $z^{\alpha}$
denote the
      monomial in noncommuting indeterminates
      $$
        z^{\alpha} = z_{i_{N}} \cdots z_{i_{1}}
      $$
      and we let $z^{\emptyset} = 1$.
      We extend this noncommutative functional calculus to a $d$-tuple of
      operators ${\mathbf A} = (A_{1}, \dots, A_{d})$ on a Hilbert space
      $\cX$:
      \begin{equation}  \label{bAv}
     {\mathbf A}^{v} = A_{i_{N}} \cdots A_{i_{1}}\quad \text{if}\quad
v = i_{N} \cdots i_{1} \in \cF_{d} \setminus \{ \emptyset\};\quad {\mathbf
      A}^{\emptyset} = I_{\cX}.
      \end{equation}

      We will also have need of the {\em
transpose operation} on $\cF_{d}$:
\begin{equation} \label{transpose}
       \alpha^{\top} = i_{1} \cdots i_{N} \quad\text{if}\quad \alpha =
i_{N}\cdots i_{1}.
\end{equation}

     A natural analogue of the Szeg\"o kernel is the noncommutative
     Szeg\"o kernel
\begin{equation}
     k_{\text{Sz}}(z,w) = \sum_{\alpha \in \cF_{d}} z^{\alpha} 
w^{\alpha^{\top}}.
\label{szego}
\end{equation}
     The associated reproducing kernel Hilbert space
     $\cH(k_{\text{Sz}})$ (in the sense of \cite{NFRKHS}) is a natural
     analogue of the classical Hardy space $H^{2}({\mathbb D})$; we
     recall all the relevant definitions and main properties more
     precisely in Section \ref{S:NFRKHS}.
     Our main purpose here is to obtain the
     analogues the Theorems \ref{T:I1}--\ref{T:I6}  above with the classical
     Szeg\"o kernel replaced by its noncommutative analogue
     \eqref{szego}.

     In particular, the analogue of Theorem  \ref{T:I6} involves
     the study of shift-invariant subspaces of the Fock space 
$H^{2}_{\cY}(\cF_{d})$
     generated by a collection of homogeneous interpolation conditions
     defined via a functional calculus with noncommutative operator
     argument.  We mention that interpolation problems in the
     noncommutative Schur-multiplier class defined by nonhomogeneous
     interpolation conditions associated with such a functional calculus
     have been studied recently by a number of authors, including the late
     Tiberiu Constantinescu to whom this paper is dedicated (see
     \cite{BB-noncomint, CJ, Popescu-Nehari, Popescu-Memoir}).
     While the Fock-space version of the Beurling-Lax theorem already
     appears in the work of Popescu \cite{PopescuNF1} (see also
     \cite{BBF1}), the proof here through inner solution of a homogeneous
     interpolation problem gives an alternative approach.

     The present paper (with the exception of the final Section
     \ref{S:BL} ) parallels our companion paper \cite{BBF2} where
     corresponding results are worked out with the noncommutative Szeg\"o
     kernel \eqref{szego}  replaced by the so-called Arveson kernel
     $k_{d}(\blam, \bzeta) = 1/(1 - \langle \blam, \bzeta
     \rangle_{{\mathbb C}^{d}})$
     which is positive on the unit ball ${\mathbb B}^{d} = \{\blam =
     (\lam_{1}, \dots, \lam_{d}) \colon \sum_{k=1}^{d} |\lam_{k}|^{2}
     < 1\}$ of ${\mathbb C}^{d}$.
      There the corresponding results are
     more delicate; in particular, the observable, coisometric
     realization for a contractive multiplier is unique only in very
     special circumstances, but the
     nonuniqueness can be explicitly characterized.  In contrast, the
     results obtained here for the setting of the noncommutative Szeg\"o
     kernel $k_{\text{Sz}}(z,w)$ parallel more directly the situation
     for the classical univariate case.

     The paper is organized as follows.  After the present Introduction,
     Section \ref{S:NFRKHS} recalls the main facts from \cite{NFRKHS}
     which are needed in the sequel.  Section \ref{S:NC-Schur}
     introduces the noncommutative Schur class of contractive Fock-space
     multipliers $S$ and the associated
     noncommutative positive kernel $K_{S}(z,w)$, and develops the
     noncommutative analogues of Theorems \ref{T:I1} and \ref{T:I2}.  In
     fact, various pieces of the noncommutative version of Theorem
     \ref{T:I1} (see theorem \ref{T:NC1} below) are already worked out
     in \cite{NFRKHS, PopescuNF2, Cuntz-scat}.  In connection with the
     noncommutative analogue of Theorem \ref{T:I2} (see Theorems
     \ref{T:CAtoS} and \ref{T:CAStoB} below), we rely on our paper
     \cite{BBF1} where the structure of noncommutative formal
     reproducing kernel spaces of the type $\cH(K_{C,A})$ were worked out.
     Section \ref{S:dBR} introduces the noncommutative functional-model
     coisometric colligation ${\mathbf U}_{dBR}$ and obtains the
     analogue of Theorem \ref{T:I3} for the Fock space setting (see
     Theorem \ref{T:NC3} below).  This functional model is the
     Brangesian model parallel to the noncommutative Sz.-Nagy-Foia\c{s}
     model for a row contraction found in \cite{PopescuNF2, Cuntz-scat}.
     The final Section \ref{S:BL} uses previous results concerning
     $\cH(K_{S})$ and $\cH(K_{C,A})$ to arrive at the Fock-space version
     of Theorem \ref{T:I6} (see Theorems \ref{T:shift=int} and
     \ref{T:BLhomint} below) in a simple way.

\section{Noncommutative formal reproducing kernel Hilbert spaces}
\label{S:NFRKHS}

We now recall some of the basic ideas from \cite{NFRKHS} concerning 
noncommutative
formal reproducing kernel Hilbert spaces.
We let $z = (z_{1}, \dots, z_{d})$, $w =
(w_{1}, \dots, w_{d})$ be two sets of noncommuting indeterminates
and we let $\cF_{d}$ be the free semigroup generated by the alphabet
$\{1, \dots, d\}$ with unit element equal to the empty word
$\emptyset$ as in the introduction. Given a coefficient Hilbert space
$\cY$ we let $\cY\langle z \rangle$ denote the set of all polynomials in
$z = (z_{1}, \dots, z_{d})$ with coefficients in $\cY:$
$$
    \cY\langle z \rangle = \left\{ p(z) = \sum_{\alpha \in \cF_{d}}
    p_{\alpha} z^{\alpha} \colon p_{\alpha} \in \cY \text{ and }
    p_{\alpha} = 0 \text{ for all but finitely many } \alpha \right\},
    $$
    while $\cY \langle \langle z \rangle \rangle$ denotes the set of all
    formal power series in the indeterminates $z$ with coefficients in
    $\cY$:
$$ \cY\langle \langle z \rangle \rangle = \left\{ f(z) = \sum_{\alpha \in
    \cF_{d}} f_{\alpha} z^{\alpha} \colon f_{\alpha} \in \cY \right\}.
    $$
Note that vectors in $\cY$ can be considered as Hilbert space
operators between ${\mathbb C}$ and $\cY$.  More generally, if $\cU$
and $\cY$ are two Hilbert spaces, we let
$\cL(\cU, \cY)\langle z \rangle$ and $\cL(\cU, \cY)\langle \langle z
\rangle \rangle$ denote the space of polynomials (respectively,
formal power series) in the noncommuting indeterminates $z = (z_{1},
\dots, z_{d})$ with coefficients in $\cL(\cU, \cY)$.
Given $S = \sum_{\alpha \in \cF_{d}} s_{\alpha}
z^{\alpha} \in \cL(\cU, \cY)\langle \langle z \rangle
\rangle$ and $f = \sum_{\beta \in \cF_{d}} f_{\beta} z^{\beta} \in
\cU\langle \langle z \rangle \rangle$, the product $S(z) \cdot f(z)
\in \cY \langle \langle z \rangle \rangle$ is defined as an element
of $\cY\langle \langle z \rangle \rangle$ via the noncommutative
convolution:
\begin{equation} \label{multiplication}
S(z) \cdot f(z) = \sum_{\alpha, \beta \in \cF_{d}} s_{\alpha}
f_{\beta} z^{\alpha \beta} =
\sum_{v \in \cF_{d}} \left( \sum_{\alpha, \beta \in \cF_{d} \colon
\alpha \cdot \beta = v} s_{\alpha} f_{\beta} \right) z^{v}.
\end{equation}
Note that the coefficient of $z^{v}$ in \eqref{multiplication} is
well defined since any given word $v \in \cF_{d}$ can be decomposed
as a product $v = \alpha \cdot \beta$ in only finitely many distinct
ways.

In general, given a coefficient Hilbert space $\cC$, we use the
$\cC$ inner product to generate a pairing
$$ \langle \cdot, \, \cdot \rangle_{\cC \times \cC\langle \langle w
\rangle \rangle} \colon \cC \times \cC\langle \langle w \rangle
\rangle \to \cC\langle \langle w \rangle \rangle
$$
via
$$
\left\langle c, \sum_{\beta \in \cF_{d}} f_{\beta} w^{\beta}
\right\rangle_{\cC \times
\cC\langle \langle w \rangle \rangle} = \sum_{\beta \in \cF_{d}} \langle
c, f_{\beta} \rangle_{\cC} w^{\beta^{\top}} \in \cC \langle \langle
w \rangle \rangle.
$$
We also may use the pairing in the reverse order
$$
\left\langle \sum_{\alpha \in \cF_{d} } f_{\alpha} w^{\alpha}, c
\right\rangle_{\cC\langle \langle w \rangle \rangle \times \cC} =
\sum_{\alpha \in \cF_{d}} \langle f_{\alpha}, c \rangle_{\cC}
w^{\alpha} \in \cC\langle \langle w \rangle \rangle.
$$
These are both special cases of the more general pairing
$$
\left\langle \sum_{\alpha \in \cF_{d}} f_{\alpha} w^{\prime \alpha},
\sum_{\beta \in \cF_{d}} g_{\beta} w^{\beta}
\right\rangle_{\cC\langle \langle w' \rangle \rangle \times \cC
\langle \langle w \rangle \rangle} =
\sum_{\alpha, \beta \in \cF_{d}} \langle f_{\alpha}, g_{\beta}
\rangle_{\cC} w^{\beta^{\top}} w^{\prime \alpha}.
$$
Suppose that $\cH$ is a Hilbert space whose elements are formal power series in
$\cY \langle \langle z \rangle \rangle$ and that $K(z,w) =
\sum_{\alpha, \beta \in \cF_{d}} K_{\alpha, \beta} z^{\alpha}
w^{\beta^{\top}}$ is a formal power series in the two sets of $d$
noncommuting indeterminates $z = (z_{1}, \dots, z_{d})$ and $w =
(w_{1}, \dots, w_{d})$.  We say that {\em $K(z,w)$ is a reproducing
kernel for $\cH$} if, for each $\beta  \in \cF_{d}$ the formal power series
$$ K_{\beta}(z): = \sum_{\alpha \in \cF_{d}} K_{\alpha, \beta}
z^{\alpha}\quad\mbox{belongs to} \; \;  \cH
$$
and we have the reproducing property
$$
     \langle f, K(\cdot, w) y \rangle_{\cH \times \cH\langle \langle w
     \rangle \rangle} = \langle f(w), y \rangle_{\cY\langle \langle w
     \rangle \rangle \times \cY}\quad\text{ for every } f\in\cH.
$$
As a consequence we then also have
$$
\langle K(\cdot,w')y', K(\cdot, w) y \rangle_{\cH\langle \langle
w'\rangle \rangle \times \cH\langle \langle w \rangle \rangle} =
\langle K(w,w') y',y\rangle_{\cY\langle \langle w,w' \rangle \rangle
\times \cY}.
$$
It is not difficult to see that a reproducing kernel for a given
$\cH$ is necessarily unique.

Let us now suppose that $\cH$ is a Hilbert space whose elements are
formal power series $f(z) = \sum_{\alpha \in \cF_{d}} f_{v} z^{v}
\in \cY \langle \langle z \rangle \rangle$ for a coefficient Hilbert
space $\cY$.  We say that $\cH$ is a NFRKHS ({\em noncommutative
formal reproducing kernel Hilbert space}) if, for each $\alpha \in
\cF_{d}$, the linear operator $\Phi_{\alpha} \colon \cH \to \cY$
defined by $f(z) = \sum_{v \in \cF_{d}} f_{v} z^{v} \mapsto
f_{\alpha}$ is continuous. In this case, define $K(z,w) \in 
\cL(\cY)\langle \langle
z, w \rangle \rangle$ by
$$
K(z,w) = \sum_{\beta \in \cF_{d}} \Phi_{\beta}^{*} w^{\beta^{\top}}
=:\sum_{\alpha, \beta \in \cF_{d}} K_{\alpha, \beta} z^{\alpha}
w^{\beta^{\top}}.
$$
Then one can check that $K(z,w)$ is a reproducing kernel for $\cH$ in
the sense defined above.  Conversely (see \cite[Theorem 3.1]{NFRKHS}),
a given formal kernel $K(z,w) = \sum_{\alpha, \beta \in \cF_{d}}
K_{\alpha, \beta} z^{\alpha} w^{\beta^{\top}}
\in \cL(\cY)\langle \langle z, w \rangle \rangle$ is the reproducing
kernel for some NFRKHS $\cH$ if and only if $K$ is positive definite
in either one of the equivalent senses:
\begin{enumerate}
       \item $K(z,w)$ has a factorization
       \begin{equation} \label{pos1}
       K(z,w) = H(z) H(w)^{*}
       \end{equation}
       for some $H \in \cL(\cX, \cY)\langle \langle z \rangle \rangle$
       for some auxiliary Hilbert space $\cX$.  Here
      $$
      H(w)^{*} = \sum_{\beta \in \cF_{d}} H_{\beta}^{*} w^{\beta^{\top}} =
      \sum_{\beta \in \cF_{d}} H_{\beta^{\top}}^{*} w^{\beta}\quad \text{if}
\quad
      H(z) = \sum_{\alpha \in \cF_{d}} H_{\alpha} z^{\alpha}.
      $$
      \item For all finitely supported $\cY$-valued functions $\alpha
      \mapsto y_{\alpha}$ it holds that
      \begin{equation} \label{pos2}
      \sum_{\alpha, \alpha' \in \cF_{d}} \langle K_{\alpha, \alpha'}
      y_{\alpha'}, y_{\alpha} \rangle \ge 0.
      \end{equation}
\end{enumerate}
If $K$ is such a positive kernel, we denote by $\cH(K)$ the
associated NFRKHS consisting of elements of $\cY\langle \langle z
\rangle \rangle$.

\section{The noncommutative Schur class: associated positive kernels
and transfer-function realization}  \label{S:NC-Schur}

A natural analogue of the vector-valued Hardy space over the unit disk
(see e.g.~\cite{PopescuNF1})
is the Fock space with coefficients in $\cY$ which we denote here by
$H^{2}_{\cY}(\cF_{d})$:
$$
    H^{2}_{\cY}(\cF_{d}) = \left\{ f(z) = \sum_{\alpha \in \cF_{d}}
    f_{\alpha} z^{v} \colon \sum_{\alpha \in \cF_{d}} \|
    f_{\alpha}\|^{2} < \infty \right\}.
$$
When $\cY = {\mathbb C}$ we write simply $H^{2}(\cF_{d})$.
As explained in \cite{NFRKHS}, $H^{2}(\cF_{d})$ is a NFRKHS with
reproducing kernel  equal to the following noncommutative analogue of the
classical Szeg\"o kernel:
\begin{equation}  \label{kSz}
       k_{\text{Sz}}(z,w) = \sum_{\alpha \in \cF_{d}} z^{\alpha}
       w^{\alpha^{\top}}.
\end{equation}
Thus we have in general 
$H^{2}_{\cY}(\cF_{d}) = \cH(k_{\text{Sz}}  \otimes I_{\cY})$.
We let $S_{j}$ denote the shift operator
\begin{equation}  \label{shift}
S_{j} \colon f(z) = \sum_{v \in \cF_{d}} f_{v} z^{v} \mapsto
f(z) \cdot z_{j} = \sum_{v \in \cF_{d}} f_{v} z^{v \cdot j} \text{
for } j = 1, \dots, d
\end{equation}
on $H^{2}_{\cY}(\cF_{d})$; when we wish to specify the coefficient
space $\cY$ explicitly, we write $S_{j} \otimes I_{\cY}$ rather than
only $S_{j}$.  The adjoint of $S_{j} \colon
H^{2}_{\cY}(\cF_{d}) \to H^{2}_{\cY}(\cF_{d})$ is then given by
\begin{equation} \label{bs}
      S_{j}^{*} \colon \sum_{v \in \cF_{d}} f_{v}z^{v} \mapsto \sum_{v \in
      \cF_{d}} f_{v \cdot j} z^{v}\quad \text{for}\quad j = 1, \dots, d.
\end{equation}

We let $\cM_{nc,d}(\cU, \cY)$ denote the set of formal power series
$S(z) = \sum_{\alpha \in \cF_{d}} s_{\alpha} z^{\alpha}$ with
coefficients $s_{\alpha} \in \cL(\cU, \cY)$ such that the associated
multiplication operator $M_{S} \colon f(z) \mapsto S(z) \cdot f(z)$
(see \eqref{multiplication}) defines a bounded operator from
$H^{2}_{\cU}(\cF_{d})$ to $H^{2}_{\cY}(\cF_{d})$. It is not difficult
to show that $\cM_{nc,d}(\cU, \cY)$ is the  intertwining space for
the two tuples ${\mathbf S} \otimes I_{\cU} = (S_{1} \otimes I_{\cU},
\dots,, S_{d} \otimes I_{\cU})$ and ${\mathbf S}\otimes I_{\cY} =
(S_{1} \otimes I_{\cY}, \dots, S_{d} \otimes I_{\cY})$:
{\em an operator $X \in \cL(\cU, \cY)$ equals $X = M_{S}$ for some $S
\in \cM_{nc,d}(\cU, \cY)$ whenever $S_{j} \otimes I_{\cY}) X = X
(S_{j}\otimes I_{\cU})$ for $j = 1, \dots, d$}
(see e.g.~\cite{PopescuNF2} where, however, the
conventions are somewhat different).
We define the
noncommutative Schur class ${\mathcal S}_{nc, d}(\cU, \cY)$ to
consist of such multipliers $S$ for which $M_{S}$ has operator norm at
most 1:
\begin{equation} \label{ncSchur}
      {\mathcal S}_{nc, d}(\cU, \cY) = \{ S \in \cL(\cU, \cY) \colon
      M_{S} \colon H^{2}_{\cY}(\cF_{d}) \to H^{2}_{\cY}(\cF_{d}) \text{
      with } \|M_{S}\|_{op} \le 1 \}.
\end{equation}
The following is the noncommutative analogue of Theorem \ref{T:I1}
for this setting.

\begin{theorem}  \label{T:NC1} Let $S(z) \in \cL(\cU, \cY) \langle
     \langle z \rangle \rangle$
    be a formal power series in $z = (z_{1}, \dots, z_{d})$ with
coefficients in $\cL(\cU, \cY)$.  Then the following are equivalent:
\begin{enumerate}
       \item $S \in \cS_{nc,d}(\cU, \cY)$, i.e., $M_{S} \colon
\cU\langle
       z\rangle \to \cY \langle \langle z \rangle \rangle$ given by
       $M_{S} \colon p(z) \to S(z) p(z)$ extends to define a
       contraction operator from $H^{2}_{\cU}(\cF_{d})$ into
       $H^{2}_{\cY}(\cF_{d})$.

       \item The kernel
       \begin{equation}  \label{KS}
     K_{S}(z, w) : = k_{\text{Sz}}(z,w) - S(z) k_{\text{Sz}}(z,w)
       S(w)^{*}
       \end{equation}
       is a noncommutative positive kernel (see \eqref{pos1} and
       \eqref{pos2}).

       \item There exists a Hilbert space $\cX$ and a unitary connection
       operator $\bU$ of the form
       \begin{equation} \label{NCcolligation}
       \bU = \begin{bmatrix} A & B \\ C & D \end{bmatrix} =
       \begin{bmatrix} A_{1} & B_{1} \\ \vdots & \vdots \\ A_{d} &
       B_{d} \\ C & D \end{bmatrix} \colon \begin{bmatrix} \cX \\ \cU
       \end{bmatrix} \to \begin{bmatrix}  \cX \\ \vdots \\ \cX \\ \cY
       \end{bmatrix}
       \end{equation}
       so that $S(z)$ can be realized as a formal power series in the
       form
       \begin{equation}   \label{NCrealization}
         S(z) = D + \sum_{j=1}^{d} \sum_{v \in \cF_{d}} C A^{v}B_{j}
         z^{v}\cdot z_{j}=
         D + C (I - Z(z) A)^{-1} Z(z) B
        \end{equation}
where we have set
\begin{equation}
Z(z)=\begin{bmatrix}z_1 I_{\cX} & \ldots &
z_dI_{\cX}\end{bmatrix},\quad
A=\begin{bmatrix} A_1 \\ \vdots \\ A_d\end{bmatrix},\quad
B=\begin{bmatrix} B_1 \\ \vdots \\ B_d\end{bmatrix}.
\label{1.6a}
\end{equation}

       \item There exists a Hilbert space $\cX$ and a contractive block
       operator matrix $\bU$
       as in \eqref{NCcolligation} such that $S(z)$ is given as in
       \eqref{NCrealization}
       \end{enumerate}
       \end{theorem}

       \begin{proof} (1) $\Longrightarrow$ (2) is Theorem 3.15 in
       \cite{NFRKHS}. A proof of (2)
       $\Longrightarrow$ (3) is done in \cite[Theorem 5.4.1]{Cuntz-scat}
       as an application of the Sz.-Nagy-Foia\c s model theory for row
       contractions worked out there following ideas of Popescu
       \cite{PopescuNF1, PopescuNF2}; an
       alternative proof via the ``lurking isometry argument'' can be
       found in \cite[Theorem 3.16]{NFRKHS}.  The implication (3)
       $\Longrightarrow$ (4) is trivial.
       The content of (4) $\Longrightarrow$ (1) amounts to Proposition
       4.1.3 in \cite{Cuntz-scat}.
       \end{proof}

       We note that formula \eqref{NCrealization} has the interpretation
       that $S(z)$ is the {\em transfer function} of the
       multidimensional linear system with evolution along $\cF_{d}$
       given by the input-state-output equations
       \begin{equation}  \label{sys}
\Sigma \colon 	\left\{ \begin{array}{ccc}
	x(1 \cdot \alpha) & = & A_{1} x(\alpha) + B_{1} u(\alpha)  \\
	\vdots &   & \vdots  \\
	x(d \cdot \alpha) & = & A_{d} x(\alpha) + B_{d} u(\alpha) \\
	y(\alpha) & = & C x(\alpha) + D u(\alpha)
	\end{array} \right.
       \end{equation}
       initialized with $x(\emptyset) = 0$.  Here $u(\alpha)$ takes
       values in the input space $\cU$, $x(\alpha)$ takes
       values in the state space $\cX$, and $y(\alpha)$ takes values in
       the output space $\cY$ for each $\alpha \in \cF_{d}$.  If we
       introduce the noncommutative $Z$-transform
       $$ \{x(\alpha)\}_{\alpha \in \cF_{d}} \mapsto \widehat{x}(z) : =
       \sum_{\alpha \in
       \cF_{d}} x(\alpha) z^{\alpha}
       $$
       and apply this transform to each of the system
       equations in \eqref{sys}, one can solve for $\widehat y(z)$ in
       terms of $\widehat u(z)$ to arrive at
       $$
         \widehat y(z) = T_{\Sigma}(z) \cdot \widehat u(z)
       $$
       where the {\em transfer function} $T_{\Sigma}(z)$ of the system
       \eqref{sys} is the formal power series with coefficients in
       $\cL(\cU, \cY)$ given by
       \begin{equation}  \label{transfunc}
        T_{\Sigma}(z) = D + \sum_{j=1}^{d} \sum_{\alpha \in \cF_{d}} C
        \bA^{v}B_{j} z^{v}  z_{j} = D + C (I - Z(z) A)^{-1} Z(z) B.
       \end{equation}
       For complete details, we refer to \cite{Cuntz-scat, BGM1, BGM2}.

       The implication (4) $\Longrightarrow$ (2) can be seen directly
       via the explicit identity \eqref{NC-ADR} given in the next
       proposition; for the commutative case we refer to \cite[Lemma
       2.2]{ADR}.

       \begin{proposition}  \label{P:NC-ADR}
     Suppose that ${\mathbf U} = \sbm{ A & B \\ C &
       D } \colon \cX \oplus \cU \to \cX^{d} \oplus \cY$ is contractive
       with associated transfer function $S \in {\mathcal S}_{nc,d}(\cU,
       \cY)$ given by \eqref{NCrealization}.  Then the kernel
       $K_{S}(z,w)$ given by \eqref{KS} is can also be represented as
\begin{equation}
K_{S}(z,w) =C(I_{\cX} - Z(z)A)^{-1} (I_{\cX} - A^{*} Z(w)^{*})^{-1} C^{*}
           + D_{S}(z,w)
\label{NC-ADR}
\end{equation}
where
\begin{align}
D_{S}(z,w) =&\begin{bmatrix} C (I_{\cX} - Z(z) A)^{-1} Z(z) & I_{\cY}
	\end{bmatrix}  k_{\text{Sz}}(z,w)\notag\\
	&\cdot (I - {\mathbf U}	{\mathbf U}^{*})
	 \begin{bmatrix} Z(w)^{*} (I - A^{*} Z(w)^{*})^{-1}
	C^{*} \\ I_{\cY} \end{bmatrix}.
	\label{DS}
     \end{align}
     \end{proposition}

     \begin{proof}
       For a fixed $\alpha \in \cF_{d}$, let us set
       \begin{align}
      	X_{\alpha}  & = z^{\alpha} w^{\alpha^{\top}} I_{\cY} - S(z)
         z^{\alpha}w^{\alpha^{\top}} S(w)^{*},  \label{xa}\\
	Y_{\alpha} & = \begin{bmatrix} C (I - Z(z)A)^{-1} Z(z) & I_{\cY}
	\end{bmatrix} z^{\alpha} w^{\alpha^{\top}} (I - {\mathbf U}
	{\mathbf U}^{*}) \begin{bmatrix} Z(w)^{*}(I - A^{*}
	Z(w)^{*})^{-1} C^{*} \\ I_{\cY} \end{bmatrix}.\notag
      \end{align}
      Note that by \eqref{KS} and \eqref{kSz},
$$
\sum_{\alpha \in \cF_{d}} X_{\alpha} = K_{S}(z,w)\quad\mbox{and}\quad
\sum_{\alpha \in \cF_{d}} Y_{\alpha}= D_{S}(z,w).
$$
Therefore \eqref{NC-ADR} is verified once we
       show that
        \begin{equation}  \label{NC-ADR'}
	 \sum_{\alpha \in \cF_{d}} X_{\alpha} - \sum_{\alpha \in
	 \cF_{d}} Y_{\alpha} =
	 C(I - Z(z) A)^{-1} (I - A^{*} Z(w)^{*})^{-1} C^{*}.
\end{equation}
Substituting \eqref{NCrealization} into \eqref{xa} gives
      \begin{align*}
          X_{\alpha}
          & = z^{\alpha} w^{\alpha^{\top}}I_{\cY} -
          [D + C (I - Z(z)A)^{-1} Z(z)B] \cdot z^{\alpha}
w^{\alpha^{\top}}
          \cdot \\
          & \qquad \cdot
          [D^{*} + B^{*} Z(w)^{*}(I - A^{*} Z(w)^{*})^{-1} C^{*}] \\
          & = z^{\alpha} w^{\alpha^{\top}}( I_{\cY} - DD^{*})
          -C(I - Z(z)A)Z(z)B D^{*}z^{\alpha} w^{\alpha^{\top}} \\
          & \qquad
         - z^{\alpha} w^{\alpha^{\top}} D B^{*} Z(w)^{*} (I -
          A^{*} Z(w)^{*})^{-1} C^{*} \\
          & \qquad - C(I - Z(z) A)^{-1} Z(z) B \cdot z^{\alpha}
w^{\alpha^{\top}}
          \cdot
          B^{*} Z(w)^{*} (I - A^{*} Z(w)^{*})^{-1} C^{*}.
       \end{align*}
       On the other hand, careful bookkeeping and use of the identity
       $$ I - {\mathbf U} {\mathbf U}^{*} =
       \begin{bmatrix} I - AA^{*} - BB^{*} & -AC^{*} - BD^{*} \\ -CA^{*}
	-DB^{*} & I - CC^{*} -DD^{*} \end{bmatrix}
       $$
       gives that
       \begin{align*}
	Y_{\alpha} & = C(I - Z(z)A)^{-1} Z(z) \cdot z^{\alpha}
	w^{\alpha^{\top}} \cdot (I - A A^{*} - BB^{*}) Z(w)^{*} (I -
	A^{*} Z(w)^{*})^{-1} C^{*} \\
	& \qquad - C(I - Z(z)A)^{-1} Z(z) (AC^{*} + B D^{*})
	z^{\alpha} w^{\alpha^{\top}}  \\
	& \qquad - z^{\alpha} w^{\alpha^{\top}} (C A^{*} + D B^{*})
	Z(w)^{*} (I - A^{*} Z(w)^{*})^{-1} C^{*} \\
	& \qquad
	+ z^{\alpha} w^{\alpha^{\top}} (I - CC^{*} -DD^{*}).
     \end{align*}
     Further careful bookkeeping then shows that
     \begin{align}
       &  X_{\alpha} - Y_{\alpha}  = z^{\alpha} w^{\alpha^{\top}} C
C^{*} + C (I -
         Z(z)A)^{-1}Z(z) A C^{*} z^{\alpha} w^{\alpha^{\top}} \notag \\
         & \qquad + z^{\alpha} w^{\alpha^{\top}} C A^{*} Z(w)^{*} (I -
         A^{*} Z(w)^{*})^{-1} C^{*} \notag \\
         & \qquad
         - C(I - Z(z) A)^{-1} Z(z) \cdot z^{\alpha} w^{\alpha^{\top}}
         \cdot (I - AA^{*}) Z(w)^{*} (I - A^{*} Z(w)^{*})^{-1} C^{*}
         \notag \\
         & = C(I - Z(z)A)^{-1} (z^{\alpha} w^{\alpha^{\top}} I_{\cX} -
         Z(z) z^{\alpha} w^{\alpha^{\top}}  Z(w)^{*})
         (I - A^{*} Z(w)^{*})^{-1} C^{*}.
         \label{Xalpha-Yalpha}
      \end{align}
      Note that
      $$ Z(z) \cdot z^{\alpha} w^{\alpha^{\top}}\cdot Z(w)^{*}
      = \sum_{k=1}^{d} z_{k} z^{\alpha} w^{\alpha^{\top}} w_{k}
      $$
      and hence
      $$
      \sum_{\alpha \in \cF_{d}\colon |\alpha| = N} Z(z)
      z^{\alpha} w^{\alpha^{\top}} Z(w)^{*}
      = \sum_{\alpha \in \cF_{d} \colon |\alpha| = N+1}
      z^{\alpha} w^{\alpha^{\top}} I_{\cX}.
      $$
Therefore,
      \begin{align}
        &  \sum_{\alpha \in \cF_{d}} z^{\alpha} w^{\alpha^{\top}} I_{\cX}
          - \sum_{\alpha \in \cF_{d}} Z(z) \notag
          z^{\alpha}w^{\alpha^{\top}} Z(w)^{*} \\
        & \qquad   =\sum_{N=0}^{\infty} \sum_{\alpha \in \cF_{d}\colon
|w| = N}
          z^{\alpha} w^{\alpha^{\top}} I_{\cX}
          - \sum_{N=1}^{\infty} \sum_{\alpha \in \cF_{d}\colon |w| = N}
          z^{\alpha}w^{\alpha^{\top}}I_{\cX}= I_{\cX}.
          \label{middle}
          \end{align}
     Summing \eqref{Xalpha-Yalpha} and combining with \eqref{middle}
     gives the result \eqref{NC-ADR'} as wanted.
        \end{proof}

       Given a $d$-tuple of operators $A_{1}, \dots, A_{d}$ on the
       Hilbert space $\cX$, we let $\bA = (A_{1}, \dots, A_{d})$ denote
       the operator $d$-tuple while
       $A$  denotes the associated column matrix as in \eqref{1.6a}
       considered as an operator from $\cX $ into $\cX^{d}$.  If $C$ is
       an operator from $\cX$ into an output space $\cY$, we say that
       $(C, \bA)$ is an output pair. The paper \cite{BBF1} studied
       output pairs and connections with the associated state-output
       noncommutative linear system \eqref{sys}.
     We are particularly interested in the case where in addition $(C,
     \bA)$ is {\em contractive}, i.e.,
\begin{equation}
 A_{1}^{*} A_{1} + \cdots + A_{d}^{*}A_{d} + C^{*} C \le I_{\cX}.
\label{cont}
\end{equation}
     In this case we have the following result.

     \begin{proposition}  \label{P:BBF1}
         Suppose that $(C, \bA)$ is a contractive
         output pair.  Then:
         \begin{enumerate}
	  \item The observability operator
         \begin{equation}  \label{ob-op}
       \cO_{C, \bA} \colon x \mapsto \sum_{\alpha \in
         \cF_{d}}( C \bA^{v} x) z^{\alpha} = C(I - Z(z) A)^{-1} x
         \end{equation}
         maps $\cX$ contractively into $H^{2}_{\cY}(\cF_{d})$.

         \item The space $\operatorname{Ran}\,  \cO_{C, \bA}$
         is a
         NFRKHS with norm given by
         $$
         \|  \cO_{C, \bA}x \|_{\cH(K_{C,\bA})}=\| Q x\|_{\cX}
         $$
         where $Q$ is the orthogonal projection onto
         $(\operatorname{Ker}\,\cO_{C,\bA})^{\perp}$ and with
         formal reproducing kernel $K_{C,A}$ given by
         \begin{equation}  \label{KCA}
         K_{C,\bA}(z,w) = C(I - Z(z)A)^{-1} (I - Z(w)^{*}A^{*})^{-1}
         C^{*}.
         \end{equation}

         \item $\cH(K_{C,\bA})$ is invariant under the backward shift
         operators $S_{j}^{*}$ given by \eqref{bs} for $j = 1, \dots, d$
         and moreover the difference-quotient inequality
         \begin{equation}  \label{DQineq}
          \sum_{j=1}^{d} \|S_{j}^{*}f\|^{2}_{\cH(K_{C, \bA})} \le
	\|f\|^{2}_{\cH(K_{C, \bA})} - \|f_{\emptyset}\|^{2}_{\cY}\quad
	\text{for all} \quad f \in \cH(K_{C, \bA})
	\end{equation}
	is satisfied.

	\item $\cH(K_{C,\bA})$ is isometrically included in
         $H^{2}_{\cY}(\cF_{d})$ if and only if in addition $\bA$ is
         {\em strongly stable}, i.e.,
         \begin{equation}  \label{stronglystable}
	  \lim_{N \to \infty} \sum_{\alpha \in \cF_{d} \colon
	  |\alpha| = N} \| \bA^{v} x \|^{2} = 0\quad \text{for all} \quad
            x\in \cX.
         \end{equation}
         \end{enumerate}
         \end{proposition}

         \begin{proof} We refer the reader to \cite[Theorem 
2.10]{BBF1} for complete
	   details of the proof.  Here we only note that the
	   backward-shift-invariance property in part (3) is a
	   consequence of the intertwining relation
	   \begin{equation}  \label{intertwine}
	       S_{j}^{*}\cO_{C, \bA} = \cO_{C,
	       \bA} A_{j} \quad\text{for} \quad j = 1, \dots, d
	    \end{equation}
    and that, in the observable case, \eqref{DQineq} is equivalent to
    the contractivity property \eqref{cont} of $(C, \bA)$.    \end{proof}

         The paper \cite{BBF1} studies the NFRKHSs $\cH(K)$ where the
         kernel $K$ has the special form $K_{C,\bA}$ for a contractive
         output pair as in \eqref{KCA}.
        Here we wish to study the noncommutative analogues of de
        Branges-Rovnyak spaces $\cH(K_{S})$ with $K_{S}$  given by
        \eqref{KS}.

        The following corollary to Proposition
        \ref{P:NC-ADR} gives a connection between kernels of the form
        $K_{C,\bA}$ for a contractive output pair $(C, \bA)$
        and kernels of the form $K_{S}$ for a noncommutative Schur-class
        multiplier $S \in {\mathcal S}_{nc, d}(\cU, \cY)$.

        \begin{corollary}  \label{C:NC-ADR}
	 Suppose that the operator $\bU$ of the form \eqref{NCcolligation}
           is contractive with
	 associated noncommutative Schur multiplier $S(z)$ given by
	 \eqref{NCrealization}.  Suppose that the associated
	 output-pair $(C,\bA)$ with $\bA = (A_{1}, \dots, A_{d})$ is
	 observable (i.e., the observability operator $\cO_{C,
	 \bA}$ given by \eqref{ob-op} is injective).
	Then the associated kernels
	 $K_{S}(z,w)$ and $K_{C,\bA}(z,w)$
	 given by \eqref{KS} and \eqref{KCA} are the same
	 \begin{equation}  \label{identical-kernels}
	   K_{S}(z,w) = K_{C,\bA}(z,w)
	 \end{equation}
	 if and only if ${\mathbf U}$ is coisometric.
     \end{corollary}

     \begin{proof}
         By Proposition \ref{P:NC-ADR} the identity of kernels
         \eqref{identical-kernels} holds if and only if the defect
         kernel $D_{S}(z,w)$ defined in \eqref{DS} is zero.
         Let us partition $I - {\mathbf U} {\mathbf
         U}^{*}$ as a $(d+1) \times (d+1)$ block matrix with respect to
         the $(d+1)$-fold decomposition $\cX^{d} \oplus \cY$ of its
         domain and range spaces
         $$
         I - {\mathbf U} {\mathbf U}^{*} = [M_{i,j}]_{1 \le i, j \le d+1}
         $$
         and let us write $D_S(z,w)$ as a formal power series
         $$
          D_S(z,w) = \sum_{v,v' \in \cF_{d}} D_{v,v'} z^{v} w^{v'}.
         $$
It follows from \eqref{DS} that $D_{v,v'}$ is given by
         \begin{align*}
	  D_{v,v'}&  = \sum_{\beta, \alpha, \gamma \in \cF_{d}, i,j \in
	  \{1, \dots, d\} \colon \beta i \alpha = v, \alpha^{\top} j
	  \gamma^{\top} = v'} C A^{\beta} M_{i,j} A^{* \gamma^{\top}}
	  C^{*} \\
	  & \qquad + \sum_{\beta \in \cF_{d}, i \in \{1, \dots, d\}
	  \colon \beta i = v (v^{\prime \top})^{-1}} M_{i,d+1} \\
	  & \qquad + \sum_{j \in \{1, \dots, d \}, \beta \in \cF_{d}
	  \colon j \gamma^{\top} = v'(v^{\top})^{-1}} M_{d+1,j}+ 
          M_{d+1,d+1},
     \end{align*}
     where in general we write
     $$ v w^{-1} = \begin{cases} v' &\text{if } v=v'w \\
               \text{undefined} &\text{otherwise.}
	      \end{cases}
    $$
     Considering the case $v=v' = \emptyset$ leads to $M_{d+1, d+1} = 0$.
     Considering next the case $v=i_{0}$, $v' = \emptyset$ leads to
     $M_{i_{0}, d+1} = 0$ for $i_{0} = 1, \dots, d$.  Similarly, the case
     $v = \emptyset$, $v' = j_{0}$ leads to $M_{d+1,j_{0}} = 0$ for
     $j_{0} =1, \dots, d$.  Considering next the case $v = i_{0}$, $v' =
     j_{0}$ leads to $C M_{i_{0},j_{0}} C^{*} = 0$ for all $i_{0},j_{0}
     = 1, \dots, d$, and hence $C (I - {\mathbf U} {\mathbf U}^{*} )
     C^{*} = 0$.  The general case together with an induction argument
     on the length of words leads to the general collapsing
     $$
       C A^{\beta}(I - {\mathbf U} {\mathbf U}^{*} )
       A^{*\gamma^{\top}}C^{*} = 0.
     $$
     The observability assumption then forces $I - {\mathbf U} {\mathbf
     U}^{*} = 0$, i.e., that ${\mathbf U}$ is coisometric as wanted.
     \end{proof}

     Alternatively, we can suppose that we know only the contractive
     output pair $(C, \bA)$ and we seek to find a noncommutative Schur
multiplier $S \in {\mathcal
     S}_{nc,d}(\cU, \cY)$ so that  \eqref{identical-kernels} holds.
     We start with a preliminary result.

     \begin{theorem} \label{T:CAtoS}
         Let $(C,\bA)$ with $C \in {\mathcal \cL}(\cX, \cY)$ be a
         contractive output-pair.  Then there exists an input space
         $\cU$ and an $S \in {\mathcal S}_{nc, d}(\cU, \cY)$ so that
         \begin{equation}  \label{KS=KCA}
         K_{S}(z,w) = K_{C, \bA}(z,w).
         \end{equation}
      \end{theorem}

         \begin{proof}  By the result of Corollary \ref{C:NC-ADR}, it
	  suffices to find an input space $\cU$ and an operator
	  $\sbm{B \\ D } \colon \cU \to \cX^{d} \oplus \cY$ so that
	  ${\mathbf U} : = \sbm{ A & B \\ C & D } \colon \cX \oplus
	  \cU \to \cX^{d} \oplus \cY$ is a coisometry.  The details
	  for such a coisometry-completion problem are carried out in
	  the proof of Theorem 2.1 in \cite{BBF2}.
     \end{proof}

     We now consider the situation where we are given a contractive
     output-pair $(C, \bA)$ and a noncommutative Schur multiplier $S \in
     {\mathcal S}_{nc, d}(\cU, \cY)$ so that \eqref{KS=KCA} holds.
\begin{lemma}\label{L:generalfact}
Let
$$
F(z) = \sum_{v \in \cF_{d}} F_{v}
z^{v}  \in {\mathcal L}(\cU, \cY)\quad \text{and}\quad
G(z) = \sum_{v \in \cF_{d}}
G_{v} z^{v} \in {\mathcal L}(\cU', \cY)
$$
be two formal power series.
Then the formal power series identity
\begin{equation}  \label{=kernels}
     F(z) F(w)^{*} = G(z) G(w)^{*}
\end{equation}
is equivalent to the existence of a (necessarily unique)
isometry $V$ from 
$$
{\mathcal D}_{V}  : = \overline{\operatorname{span}}_{v \in
      \cF_{d}} \operatorname{Ran}\, F_{v}^{*} \subset 
\cU\quad\mbox{onto}\quad
{\mathcal R}_{V}  : = \overline{\operatorname{span}}_{v \in
      \cF_{d}} \operatorname{Ran}\, G_{v}^{*} \subset \cU'
$$
so that the identity of formal power series
\begin{equation}  \label{series-id}
       V F(w)^{*} = G(w)^{*}
    \end{equation}
    holds.
\end{lemma}
\begin{proof} If there is an isometry $V$ satisfying
      \eqref{series-id}, equating coefficients of $v^{\top}$ gives
$$ V F_{v}^{*} = G_{v}^{*}.
$$
The isometric property of $V$ then leads to
\begin{equation} \label{coef-id}
     F_{v'} F_{v}^{*} = G_{v'} G_{v}^{*}\quad \text{for all} \quad v,v' \in
\cF_{d}
\end{equation}
from which we get
$$ \sum_{v',v \in \cF_{d}} F_{v'} F_{v}^{*} z^{v'} w^{v^{\top}} =
\sum_{v',v \in \cF_{d}} G_{v'} G_{v}^{*} z^{v'} w^{v^{\top}}
$$
which is the same as \eqref{=kernels} written out in coefficient form.

Conversely, the assumption \eqref{=kernels} leads to \eqref{coef-id}.
Then the formula
\begin{equation}  \label{formula}
     V \colon F_{v}^{*} y \mapsto G_{v}^{*} y \quad\text{for}\quad v \in
\cF_{d} \;    \text{ and } \; y \in \cY
\end{equation}
extends by linearity and continuity to a well-defined isometry (still
denoted by $V$) from ${\mathcal D}_{V}$ onto ${\mathcal R}_{V}$.
Since identification of coefficients of $z^{v}$ on both sides of
\eqref{series-id} reduces to \eqref{formula}, we see that
\eqref{series-id} follows as wanted.
\end{proof}

     \begin{lemma}  \label{L:2.1}
         Let $(C,\bA)$ be a contractive output pair and $S \in {\mathcal
         L}(\cU, \cY)\langle \langle z \rangle \rangle$ a formal power
         series.  Then the following are equivalent:
         \begin{enumerate}
	  \item \eqref{KS=KCA} holds, i.e.,
	\begin{equation}  \label{KS=KCA'}
	    C(I - Z(z)A)^{-1} (I - A^{*} Z(w)^{*})^{-1} C^{*} =
	    k_{\text{Sz}}(z,w) I_{\cY} -
	    S(z) k_{\text{Sz}}(z,w) S(w)^{*}.
	  \end{equation}

	  \item The alternative version of \eqref{KS=KCA'} holds:
	  \begin{equation}  \label{KS=KCAalt}
	      C(I - Z(z)A)^{-1} (I_{\cX} - Z(z) Z(w)^{*})  (I - A^{*}
	      Z(w)^{*})^{-1} C^{*} = I - S(z) S(w)^{*}.
	   \end{equation}

	   \item There is an isometry
$$
V = \begin{bmatrix}A_{V} & B_{V} \\
	   C_{V} & D_{V}\end{bmatrix} \colon
[\overline{\operatorname{Ran}}
	   (\cO_{C,\bA})^{*}]^{d} \oplus \cY \to \cX \oplus\cU
$$
   so that we have the identity of formal power series:
	   \begin{equation}  \label{isomV}
	       \begin{bmatrix}  A_{V} & B_{V} \\ C_{V} & D_{V}
	       \end{bmatrix}  \begin{bmatrix} Z(w)^{*} (I -
	       A^{*} Z(w)^{*} C^{*} \\ I_{\cY} \end{bmatrix} =
	       \begin{bmatrix} (I - A^{*} Z(w)^{*})^{-1} C^{*} \\
		   S(w)^{*} \end{bmatrix}.
	  \end{equation}
	  \end{enumerate}
\end{lemma}

\begin{proof}
       \textbf{(1) $\Longleftrightarrow$ (2):}  Suppose that
\eqref{KS=KCA'} holds. Then
       \begin{align*}
&	C(I - Z(z) A)^{-1} Z(z) Z(w)^{*} (I - A^{*} Z(w)^{*})^{-1}
	C^{*} \\
	& \qquad =
	\sum_{k=1}^{d}  w_{k} C (I - Z(z)A)^{-1} (I -
	A^{*}Z(w)^{*})^{-1} C^{*} z_{k} \\
    & \qquad = \sum_{k=1}^{d} w_{k} k_{\text{Sz}}(z,w) z_{k} - S(z)
\left(
    \sum_{k=1}^{d} w_{k} k_{\text{Sz}}(z,w) z_{k} \right) S(w)^{*} \\
    & \qquad = (k_{\text{Sz}}(z,w) - 1)I_{\cY} - S(z) \left(
k_{\text{Sz}}(z,w)-1\right) S(w)^{*}
    \end{align*}
and consequently,
    \begin{align*}
        & C(I - Z(z)A)^{-1} (I - Z(z) Z(w)^{*}) (I -
        A^{*}Z(w)^{*})^{-1}C^{*} \\
        & \qquad =
        k_{\text{Sz}}(z,w) I_{\cY} - S(z)k_{\text{Sz}}(z,w) S(w)^{*} \\
    & \qquad \qquad   - \left[  (k_{\text{Sz}}(z,w) - 1) I_{\cY} - S(z)
        \left(k_{\text{Sz}}(z,w)-1\right) S(w)^{*} \right] \\
     &\qquad = I_{\cY} - S(z) S(w)^{*}
     \end{align*}
     and we recover \eqref{KS=KCAalt} as desired.

     Conversely, assume that \eqref{KS=KCAalt} holds.  Multiplication
     of \eqref{KS=KCAalt} on the left by $w^{\gamma^{\top}}$ and on the
     right by $z^{\gamma}$ gives
     \begin{align}
         & C(I - Z(z)A)^{-1} \left(z^{\gamma}w^{\gamma^{\top}}I_{\cX} -
Z(z)z^{\gamma}w^{\gamma^{\top}}
         Z(w)^{*}\right) (I - A^{*}Z(w)^{*})^{-1}C^{*} \notag \\
           & \qquad = z^{\gamma} w^{\gamma^{\top}}I_{\cY} - S(z)
         z^{\gamma}w^{\gamma^{\top}} S(w)^{*}.
         \label{KS-KCApre}
     \end{align}
     Summing up \eqref{KS-KCApre} over all $\gamma \in \cF_{d}$ leaves
     us with \eqref{KS=KCA'}.  This completes the proof of (1)
     $\Longleftrightarrow$ (2).

     \textbf{(2) $\Longleftrightarrow$ (3):}
     Observe that \eqref{KS=KCAalt} can be written in equivalent block matrix
     form as
     \begin{align*}
         & \begin{bmatrix} C(I - Z(z)A)^{-1} Z(z) & I_{\cY} \end{bmatrix}
         \begin{bmatrix} Z(w)^{*}(I - A^{*}Z(w)^{*})^{-1}C^{*} \\
	  I_{\cY} \end{bmatrix} \\
	  & \qquad =
	  \begin{bmatrix} C (I - Z(z) A)^{-1} & S(z) \end{bmatrix}
	      \begin{bmatrix} (I - A^{*} Z(w)^{*})^{-1} C^{*}  \\
		  S(w)^{*} \end{bmatrix}.
\end{align*}
Now we apply Lemma \ref{L:generalfact} to the particular case
$$
F(w)^{*} = \begin{bmatrix} Z(w)^{*}(I - A^{*}Z(w)^{*})^{-1}C^{*} \\
	  I_{\cY} \end{bmatrix}, \qquad
G(w)^{*}  = \begin{bmatrix} (I - A^{*} Z(w)^{*})^{-1} C^{*}  \\
		  S(w)^{*} \end{bmatrix}
$$
to see the equivalence of  (2) and (3).   It is easily
checked that ${\mathcal D}_{V}$ for our case here is the $d$-fold
inflation of the observability subspace inside $\cX^{d}$:
$$
{\mathcal D}_{V} =
[\overline{\operatorname{span}}_{v \in \cF_{d}} \operatorname{Ran}\,
\bA^{* v} C^{*}]^{d} \oplus \cY = [ \overline{\operatorname{Ran}} \,
\cO_{C, \bA})^{*}]^{d} \oplus \cY.
$$
\end{proof}

       \begin{theorem}  \label{T:CAStoB}
	Suppose that $S(z) \in \cS_{nc, d}(\cU, \cY)$ and that $(C,
	\bA)$ is an observable,
	contractive output-pair  such that \eqref{KS=KCA} holds.
     Then there exists a unique operator
     $ B = \sbm{ B_{1} \\ \vdots \\ B_{d}} \colon \cU
     \to \cX^{d}$
     so that  $\bU = \sbm{ A & B \\ C & s_{\emptyset}}$ is  a coisometry
and
     $\bU$ provides a realization for $S$: $S(z) =  s_{\emptyset} + C
     (I - Z(z) A)^{-1} Z(z) B$.
     \end{theorem}

     \begin{proof} We are given the operators $A \colon \cX \to
         \cX^{d}$, $C \colon \cX \to \cY$ and $D = S_{\emptyset} \colon
         \cU \to \cY$ and seek an operator $B \colon \cU \to \cX^{d}$ so
         that $S(z) = D + C(I - Z(z)A)^{-1} Z(z)B$, or, what is the same,
         so that
         $$
         S(w)^{*} = D^{*} + B^{*}Z(w)^{*}(I - A^{*} Z(w)^{*})^{-1} C^{*}.
         $$
      This last identity can be rewritten as
      \begin{equation}  \label{want1}
          \begin{bmatrix} A^{*} & C^{*} \\ B^{*} & D^{*} \end{bmatrix}
	   \begin{bmatrix} Z(w)^{*}(I - A^{*} Z(w)^{*})^{-1}C^{*} \\
	       I_{\cY} \end{bmatrix} = \begin{bmatrix} (I - A^{*}
	       Z(w)^{*})^{-1} C^{*} \\ S(w)^{*} \end{bmatrix}
     \end{equation}
     since the identity
     $$
     A^{*} Z(w)^{*}(I - A^{*} Z(w)^{*})^{-1}C^{*} + C^{*} = (I - A^{*}
     Z(w)^{*})^{-1} C^{*}
     $$
     expressing equality of the top components holds true automatically.
     Lemma \ref{L:2.1} tells us that there is an isometry $V =
     \sbm{A_{V} & B_{V} \\ C_{V} & D_{V}} \colon
      \cX^{d} \oplus  \cY \to
     \cX \oplus \cU$ which has the same action as desired by $\sbm{A^{*}
     & C^{*} \\ B^{*} & D^{*} }$ in \eqref{want1}.  It suffices to set
     $B^{*} = C_{V}$.
     \end{proof}

     We say that two colligations ${\mathbf U} = \sbm{A & B \\ C & D } \colon
     \cX \oplus \cU \to \cX^{d} \oplus \cY$ and ${\mathbf U}' = \sbm{A'
     & B' \\ C' & D' }$ are {\em unitarily equivalent} if there is a
     unitary operator $U \colon \cX \to \cX'$ such that
     $$
     \begin{bmatrix} \oplus_{k=1}^{d} U & 0 \\ 0 & I_{\cY} \end{bmatrix}
         \begin{bmatrix} A & B \\ C & D \end{bmatrix} = \begin{bmatrix}
	  A' & B' \\ C' & D' \end{bmatrix} \begin{bmatrix} U & 0 \\ 0
	  & I_{\cU} \end{bmatrix}.
     $$

     \begin{corollary} \label{C:unique-coisometric}
         Any two observable, coisometric realizations ${\mathbf U}$ and
         ${\mathbf U}'$ for the same $S \in \cS_{nc,d}(\cU, \cY)$ are
         unitarily equivalent.
      \end{corollary}

      \begin{proof}
          Suppose that ${\mathbf U} = \sbm{A & B \\ C & D }$ and
          ${\mathbf U}' = \sbm{A' & B' \\ C' & D'}$ are two such
          realizations.  From Proposition \ref{P:NC-ADR} we see that
          $$
           K_{C, \bA}(z,w) = K_{C',\bA'}(z,w).
      $$
      Then Theorem 2.13 of \cite{BBF1} implies that $(C,\bA)$ is
      unitarily equivalent to $(C', \bA')$, so there is a unitary
      operator $U \colon \cX \to \cX'$ such that
      $$
      C' = C U^{*}\quad\text{and}\quad A_{j}' = U A_{j} U^{*}\quad
\text{for}\quad j=1,
      \dots, d.
      $$
      Then $\widetilde{\mathbf U} = \sbm{A' & (\oplus_{k=1}^{d} U) B
      \\ C' & D }$ and ${\mathbf U}' = \sbm{A' & B' \\ C' & D }$ both
      give coisometric realizations of $S$ with the same observable
      output pair $(C', \bA')$.  By the uniqueness assertion of Theorem
      \ref{T:CAtoS}, it follows that $B' = (\oplus_{k=1}^{d}U) B$ as
      well, and hence ${\mathbf U}$ and ${\mathbf U}'$ are unitarily
      equivalent.
      \end{proof}

     \section{de Branges-Rovnyak model colligations} \label{S:dBR}

     In this section we show that any $S \in \cS_{nc,d}(\cU, \cY)$ has a
     canonical observable, coisometric realization which uses
     $\cH(K_{S})$ as the state space.  We first need some preliminaries
     concerning the finer structure of the
     noncommutative de Branges-Rovnyak functional-model spaces
     $\cH(K_{S})$.  Let us denote the Taylor coefficients of
     $S(z)$ as $s_{v}$, so
	       $$
		 S(z) = \sum_{v \in \cF_{d}} s_{v} z^{v},
	       $$
to avoid confusion with the (right) shift operators $S_{j}
\colon f(z) \mapsto f(z) \cdot z_{j}$.

     Just as in the classical case, the de Branges-Rovnyak space
     $\cH(K_{S})$ has several equivalent characterizations.

     \begin{proposition}  \label{P:HKS-char}
         Let $S \in \cS_{nc, d}(\cU, \cY)$ and let $\cH$ be a Hilbert
         space of formal power series in $\cY\langle \langle z \rangle
         \rangle$.  Then the following are equivalent.

         \begin{enumerate}
	   \item $\cH$ is equal to the NFRKHS $\cH(K_{S})$
	   isometrically, where $K_S(z,w)$ is the noncommutative
	   positive kernel given by \eqref{KS}.

	   \item $\cH = \operatorname{Ran}\, (I - M_{S}
	   M_{S}^{*})^{1/2}$ with lifted norm
	  \begin{equation}  \label{lifted-norm}
	  \| (I - M_{S}M_{S})^{1/2} g \|_{\cH} = \| Q g
	  \|_{H^{2}_{\cY}(\cF_{d})}
	  \end{equation}
	  where $Q$ is the orthogonal projection of
	  $H^{2}_{\cY}(\cF_{d})$ onto $(\operatorname{Ker}\, (I -
	  M_{S}M_{S}^{*})^{1/2})^{\perp}$.

	  \item $\cH$ is the space of all formal power series $f(z)
	  \in \cY\langle \langle z \rangle \rangle$ with finite
	  $\cH$-norm, where the $\cH$-norm is given by
	\begin{equation}  \label{cH-norm}
	    \|f\|^{2}_{\cH} = \sup_{g \in H^{2}_{\cU}(\cF_{d})}
	    \left\{
	    \| f + M_{S} g \|^{2}_{H^{2}_{\cY}(\cF_{d})} - \| g
	    \|^{2}_{H^{2}_{\cU}(\cF_{d})} \right\}.
	 \end{equation}
	 \end{enumerate}
	 \end{proposition}

	 \begin{proof}
	     \textbf{ (1) $\Longleftrightarrow$ (2):}
    It is straightforward to verify the identity
$$
(I - M_{S} M_{S}^{*}) (k_{\text{Sz}}(\cdot, w) y)=
K_{S}(\cdot, w) y \text{ for each } y \in \cY.
$$
(The interpretation for this is that, for each word
$\gamma$, the coefficient of $w^{\gamma}$ of the left
hand side agrees with the coefficient of $w^{\gamma}$ on
the right hand side as elements of $\cH(k_{\text{Sz}}
I_{\cY}) = H^{2}_{\cY}(\cF_{d})$---see \cite{NFRKHS}).
We then see that
\begin{align*}
   &   \langle (I - M_{S}M_{S}^{*}) k_{\text{Sz}}( \cdot, w') y', (I -
      M_{S} M_{S}^{*}) k_{\text{Sz}}(\cdot, w) y
      \rangle_{\cH(K_{S})\langle \langle w' \rangle \rangle
      \times \cH(K_{S})\langle \langle w \rangle \rangle)} \\
   & \qquad = \langle K_{S}(\cdot, w') y', K_{S}(\cdot, w)y
   \rangle_{\cH(K_{S})\langle \langle w' \rangle \rangle
   \times \cH(K_{S}) \langle \langle w \rangle \rangle} \\
   & \qquad = \langle K_{S}(w,w')y', y \rangle_{\cY\langle \langle w',
   w \rangle \rangle \times \cY} \\
   & \qquad = \langle K_{S}(\cdot, w') y', k_{\text{Sz}}(\cdot, w)y
   \rangle_{H^{2}_{\cY}(\cF_{d})\langle \langle w' \rangle \rangle
   \times H^{2}_{\cY}(\cF_{d})\langle \langle w \rangle \rangle} \\
   & \qquad =
   \langle (I-M_{S}M_{S}^{*}) k_{\text{Sz}}(\cdot, w') y',
   k_{\text{Sz}}(\cdot, w) y \rangle_{(H^{2}_{\cY}(\cF_{d})
\langle \langle w' \rangle \rangle \times H^{2}_{\cY}(\cF_{d})
   \langle \langle w \rangle \rangle)}.
   \end{align*}
   It follows that $\operatorname{Ran}\, (I - M_{S} M_{S}^{*}) \subset
   \cH(K_{S})$ with
   $$
    \langle (I - M_{S}M_{S}^{*}) g, (I - M_{S} M_{S}^{*}) g' \rangle
    _{\cH(K_{S})} = \langle (I - M_{S}M_{S}^{*}) g, g'
    \rangle_{H^{2}_{\cY}(\cF_{d})}
   $$
   for $g, g' \in H^{2}_{\cY}(\cF_{d})$.  The precise characterization
   $\cH(K_{S}) = \operatorname{Ran} \, (I - M_{S} M_{S}^{*})^{1/2}$
   with the lifted norm \eqref{lifted-norm} now follows via a
   completion argument.

   \textbf{(2) $\Longleftrightarrow$ (3):}  This follows from the
   argument in \cite[NI-6]{sarasonbook}.

   \end{proof}

     \begin{proposition} \label{P:H(KS)} Suppose that $S \in {\mathcal
S}_{nc, d}(\cU, \cY)$
         and let $\cH(K_{S})$ be the associated NFRKHS where $K_{S}$ is
         given by \eqref{KS}.  Then the following conditions hold:
         \begin{enumerate}
	  \item
	  The NFRKHS $\cH(K_{S})$ is contained
	  contractively in $H^{2}_{\cY}(\cF_{d})$:
	  $$ \| f \|^{2}_{\cH_{\cY}(\cF_{d})} \le
	  \|f\|^{2}_{\cH(K_{S})}\quad \text{ for all }  f \in \cH(K_{S}).
	  $$

	  \item $\cH(K_{S})$ is invariant under each of the
	  backward-shift operators $S_{j}^{*}$
	  given by \eqref{bs} for $j = 1, \dots d$, and moreover,
	  the difference-quotient
	  inequality \eqref{DQineq} holds for $\cH(K_{S})$:
	  \begin{equation}  \label{DQineq-HKS}
	      \sum_{j=1}^{d} \| S_{j}^{*}f\|^{2}_{\cH(K_{S})}
	      \le \|f\|^{2}_{\cH(K_{S})} - \|f_{\emptyset}\|^{2}.
	  \end{equation}

	  \item For each $u \in \cU$ and $j = 1, \dots, d$, the
	  vector $S_{j}^{*} (M_{S}u)$ belongs to $\cH(K_{S})$ with the
	  estimate
	  \begin{equation} \label{estimate}
	      \sum_{j=1}^{d} \| S_{j}^{*} (M_{S}u) \|^{2}_{\cH(K_{S})}
	       \le \| u \|^{2}_{\cU} - \| s_{\emptyset} u \|^{2}_{\cY}.
	   \end{equation}
	  \end{enumerate}
	  \end{proposition}

	  \begin{proof}  We know from Theorem \ref{T:NC1} that $S(z)$
	      can be realized as in \eqref{NCcolligation} and
	      \eqref{NCrealization} with ${\mathbf U} = \sbm{A & B \\
	      C & D }$ a coisometry (or even unitary).
	      From Proposition \ref{P:NC-ADR} it follows that
	      $K_{S}(z,w) = K_{C, \bA}(z,w)$ and hence $\cH(K_{S}) =
	      \cH(K_{C,\bA})$ isometrically.  Conditions (1) and (2)
	      now follow from the properties of $\cH(K_{C, \bA})$ listed
	      in Proposition \ref{P:BBF1} and the discussion
	      immediately following.

	     One can also prove points (1) and (2)  directly
	     from the characterization of $\cH(K_{S})$ in part (3) of
	     Proposition \ref{P:HKS-char} (and thereby bypass
	     realization theory) as follows; these proofs
	     follow the proofs for the classical case in \cite{dbr1,
	     dbr2}.  For the contractive inclusion property (part (1)),
	     note that
	     \begin{align*}
	     \|f\|^{2}_{H^{2}_{\cY}(\cF_{d})}& =
	     \left. \left[ \|f + M_{S}g \|^{2}_{H^{2}_{\cY}(\cF_{d})} - \|
	     g \|^{2}_{H^{2}_{\cU}(\cF_{d})} \right] \right|_{g = 0}
	     \\
	     & \qquad \le
	     \sup_{g \in H^{2}_{\cU}(\cF_{d})} \left\{ \| f +
	     M_{S}g\|^{2}_{H^{2}_{\cY}(\cF_{d})} - \| g
	     \|^{2}_{H^{2}_{\cU}(\cF_{d})}\right\}
	     = \| f\|^{2}_{\cH(K_{S})}.
	     \end{align*}

	     To verify part (2), we compute
	     \begin{align*}
	&	 \sum_{j=1}^{d} \| S_{j}^{*}f \|^{2}_{\cH(K_{S})}  =
		 \sup_{g_{j}} \left\{
		 \sum_{j=1}^{d} \left[\|S_{j}^{*}f + M_{S} g_{j}
		 \|^{2}_{H^{2}_{\cY}(\cF_{d})} - \| g_{j}
		 \|^{2}_{H^{2}_{\cU}(\cF_{d})} \right]\right \} \\
		 & = \sup_{g_{j}}
		 \left\{ \sum_{j=1}^{d}\left[ \| S_{j}S_{j}^{*}f +
		 M_{S} (g_{j}z_{j}) \|^{2} - \| g_{j}z_{j}
		 \|^{2}_{H^{2}_{\cU}(\cF_{d})} \right]\right\} \\
		 & =  \sup_{g_{j}}
		 \left\{\sum_{j=1}^{d} \| S_{j}
		 S_{j}^{*}f + M_{S} (g_{j}z_{j})
		 \|^{2}_{H^{2}_{\cY}(\cF_{d})} +
		 \|f_{\emptyset}\|^{2}_{\cY} - \sum_{j=1}^{d} \|g_{j}
		 z_{j}\|^{2}_{H^{2}_{\cU}(\cF_{d})} \right\} -
		 \| f_{\emptyset}\|^{2}_{\cY} \\
		 & = \sup_{g \in H^{2}_{\cY}(\cF_{d}) \colon
		 g_{\emptyset} = 0} \left\{ \| f + M_{S} g
		 \|^{2}_{H^{2}_{\cY}(\cF_{d})} -
		 \| g \|^{2}_{H^{2}_{\cU}(\cF_{d})}\right\} -
		 \|f_{\emptyset}\|^{2}_{\cY} \\
		 & \le \sup_{g \in H^{2}_{\cY}(\cF_{d})}
		 \left\{ \| f + M_{S} g
		 \|^{2}_{H^{2}_{\cY}(\cF_{d})} -
		 \| g \|^{2}_{H^{2}_{\cU}(\cF_{d})}\right\} -
		 \|f_{\emptyset}\|^{2}_{\cY}
		 = \| f \|^{2}_{\cH(K_{S})} - \|
		 f_{\emptyset}\|^{2}_{\cY}
	\end{align*}
	and part (2) of Proposition \ref{P:H(KS)} follows.

To verify part (3), we again use the third characterization
of $\cH(K_{S})$ in Proposition \ref{P:HKS-char}. Pick
$g_1,\ldots,g_d\in H^{2}_{\cU}(\cF_{d})$ and let
	$$
	\widetilde{g} = \sum_{j=1}^{d} g_{j} z_{j}=\sum_{j=1}^{d}S_j g_{j}.
	$$
	Since $S_j^*S_i=\delta_{ij}I$ for $i,j=1,\ldots,d$ where $\delta_{ij}$
	is the Kronecker's symbol, we have
	\begin{equation}
	\|\widetilde{g}\|^2_{H^{2}_{\cU}(\cF_{d})}=\sum_{j=1}^d
	\|S_jg_j\|^2_{H^{2}_{\cU}(\cF_{d})}=\sum_{j=1}^d
	\|g_j\|^2_{H^{2}_{\cU}(\cF_{d})}
	\label{july1}
	\end{equation}
	and, since the multiplication operator $M_S$ commutes with $S_j$ for
	$j=1,\ldots,d$, we have also
	\begin{equation}
	\|M_S\widetilde{g}\|^2_{H^{2}_{\cU}(\cF_{d})}=
	\sum_{j=1}^d\|M_Sg_j\|^2_{H^{2}_{\cU}(\cF_{d})}.
	\label{july2}
	\end{equation}
	Next we note that
	\begin{align*}
        & \| S_{j}^{*}(M_{S}u) + M_{S} g_{j}\|^{2}_{H^{2}_{\cY}(\cF_{d})} =
	\|S_j^*(M_Su)\|^2+2\Re \langle S_j^*(M_Su), \, M_Sg_j\rangle+
	\|M_Sg_j\|^2 \\
   & \qquad \qquad =\langle S_jS_j^*(M_{S}u), \, M_{S}u\rangle+2\Re 
\langle M_Su, \,
	M_Sg_jz_j\rangle+\|M_Sg_j\|^2.
	\end{align*}
	Summing up the latter equalities for $j=1,\ldots,d$, making use of
	\eqref{july2}  and applying the identity
	$$
	f- f_{\emptyset}=\sum_{j=1}^d S_jS_j^*f \qquad (f\in 
H^{2}_{\cY}(\cF_{d})
	$$
	to $f=M_Su$, we get
	\begin{align}
	\sum_{j=1}^{d}\| S_{j}^{*}(M_{S}u) + M_{S}
	g_{j}\|^{2}_{H^{2}_{\cY}(\cF_{d})}&=
	\langle M_Su-s_{\emptyset}u, \, M_Su\rangle+2\Re \langle M_Su, \,
	M_S\widetilde{g}\rangle+\|M_S\widetilde{g}\|^2\notag \\
	&=\| M_Su\|^2-\|s_{\emptyset}u\|^2+2\Re \langle M_Su, \,
	M_S\widetilde{g}\rangle+\|M_S\widetilde{g}\|^2\notag \\
	&=\|M_Su+M_S\widetilde{g}\|^2_{H^{2}_{\cY}(\cF_{d})}-
	\|s_{\emptyset}u\|^2_{\cY}.\label{july3}
	\end{align}
	Since $\| M_{S}\|_{op} \le 1$ and since $\widetilde{g}_{\emptyset}=0$,
	we have
	\begin{align}
	\|M_Su+M_S\widetilde{g}\|^2_{H^{2}_{\cY}(\cF_{d})}&=
	\|M_S(u+\widetilde{g})\|^2_{H^{2}_{\cY}(\cF_{d})}\notag\\
	&\le
	\|u+\widetilde{g}\|^2_{H^{2}_{\cU}(\cF_{d})}=
	\|u\|^2_{\cU}+\|\widetilde{g}\|^2_{H^{2}_{\cU}(\cF_{d})}.\label{july4}
	\end {align}
	Adding \eqref{july1}, \eqref{july3} and \eqref{july4} gives
	$$
	\sum_{j=1}^{d}\left[ \| S_{j}^{*}(M_{S}u) + M_{S}
		     g_{j}\|^{2}_{H^{2}_{\cY}(\cF_{d})} - \|
		     g_{j}\|^{2}_{H^{2}_{\cU}(\cF_{d})} \right]
	\le \|u\|^2_{\cU}-\|s_{\emptyset}u\|^2_{\cY}.
	$$
	The latter estimate is uniform with respect to $g_j$'s and then taking
	suprema we conclude (by the third characterization of $\cH(K_{S})$ in
	Proposition \ref{P:HKS-char}) that $S_{j}^{*}(M_{S}u) \in 
\cH(K_{S})$ for
	each $j  = 1, \dots, d$ with the estimate
	$$
	\sum_{j=1}^{d} \| S_{j}^{*} (M_{S}u) \|^{2}_{\cH(K_{S})}
		       \le \| u \|^{2}_{\cU}-\|s_{\emptyset}u\|^2_{\cY},
	$$
	This concludes the proof of  Proposition \ref{P:H(KS)}.
\end{proof}

Let us define an operator $E \colon H^{2}_{\cY}(\cF_{d}) \to \cY$ by
   \begin{equation}  \label{defE}
   E \colon \sum_{v \in \cF_{d}} f_{v} z^{v} \mapsto f_{\emptyset}.
   \end{equation}
As is observed in \cite[Proposition 2.9]{BBF1} and can be
observed directly,
\begin{equation}  \label{model-obs}
E {\mathbf S}^{*v} f = E \left( \sum_{\alpha \in \cF_{d}} f_{\alpha
v^{\top}} z^{\alpha} \right) = f_{v^{\top}} \text{ for all }
f(z) = \sum_{\alpha \in \cF_{d}} f_{\alpha} z^{\alpha} \in
\cH^{2}_{\cY}(\cF_{d}) \text{ and } v \in \cF_{d}.
\end{equation}
Hence the observability operator $\cO_{E, {\mathbf S}} \colon
H^{2}_{\cY}(\cF_{d}) \to H^{2}_{\cY}(\cF_{d})$ defined as in
\eqref{ob-op} works out to be
$$ \cO_{E, {\mathbf S}^{*}} = \tau
$$
where $\tau$ is the involution on $H^{2}_{\cY}(\cF_{d})$ given by
\begin{equation}   \label{tau}
	 \tau \colon \sum_{v \in \cF_{d}} f_{v} z^{v} \mapsto
	 \sum_{v \in \cF_{d}} f_{v^{\top}} z^{v}.
	 \end{equation}
For this reason we use the ``reflected'' de Branges-Rovnyak space
	 \begin{equation}
	 \DB = \tau \circ \cH(K_{S}): = \{ \tau(f) \colon f \in \cH(K_{S})\}
	 \label{july5}
	 \end{equation}
	 as the state space for our
	 de Branges-Rovnyak-model realization of $S$ rather than simply
	 $\cH(K_{S})$ as in the classical case.  We define
	 $$
	 \| \tau(f) \|_{\DB} = \| f \|_{\cH(K_{S})}.
	 $$
	 Recall that the operator of multiplication on the right by 
the variable
	 $z_{j}$ on $H^{2}_{\cY}(\cF_{d})$ was denoted in 
\eqref{shift} by $S_{j}$
	 rather than by $S^{R}_{j}$ for simplicity.  We shall now need
	 its left counterpart, denoted by $S^{L}_{j}$ and given by
	 \begin{equation}  \label{SL}
	      S^{L}_{j} \colon f(z) = \sum_{v \in \cF_{d}} f_{v} z^{v} \mapsto
	      z_{j} \cdot f(z) = \sum_{v \in \cF_{d}} f_{v} z^{j \cdot v}
	 \end{equation}
	 with adjoint (as an operator on $H^{2}_{\cY}(\cF_{d})$) given by
	 \begin{equation} \label{SL-adj}
	      (S^{L}_{j})^{*} \colon \sum_{v \in \cF_{d}} f_{v} z^{v} \mapsto
	      \sum_{v \in \cF_{d}} f_{j \cdot v} z^{v}.
	 \end{equation}
	 For emphasis we now write $S_{j}^{R}$ rather than simply $S_{j}$.
	 We then have the following result.

	     \begin{theorem} \label{T:NC3}
		 Let $S(z) \in \cS_{nc,d}(\cU, \cY)$ and let $\DB$
	   be the associated de Bran\-ges-Rovnyak space given by
	 \eqref{july5}. Define operators
	     \begin{align*}
	     & A_{\text{dBR},j} \colon \; \DB \to \DB,\quad
	    & B_{\text{dBR},j}& \colon \; \cU \to \DB \quad (j =
	     1, \dots, d), \\
	    &  C_{\text{dBR}} \colon \; \DB \to \cY, \qquad
	    & D_{\text{dBR}}& \colon \; \cU \to \cY&
	     \end{align*}
	     by
	     \begin{align}
		& A_{\text{dBR},j} = (S^{L}_{j})^{*}|_{\DB} ,
	 \quad
		& B_{\text{dBR},j} &=  \tau  (S^{R}_{j})^{*} M_{S}|_{\cU}
		    = (S^{L}_{j})^{*} \tau M_{S}|_{\cU},
	     \notag \\
		&  C_{\text{dBR}} = E|_{\DB}, \quad
		& D_{\text{dBR}} &= s_{\emptyset}
		 \label{dBRops}
	      \end{align}
	      where $E$ is given by \eqref{defE},
		 and set
		 $$ A_{\text{dBR}} = \begin{bmatrix} A_{\text{dBR},1} \\ \vdots
		 \\ A_{\text{dBR},d} \end{bmatrix} \colon \DB \to
	 \DB^{d}, \quad
		 B_{\text{dBR}} = \begin{bmatrix} B_{\text{dBR},1} \\ \vdots \\
		 B_{\text{dBR},d} \end{bmatrix} \cU \to \DB^{d}.
	      $$
	      Then
	      $$ {\mathbf U}_{dBR} = \begin{bmatrix} A_{\text{dBR}} &
	 B_{\text{dBR}} \\
	      C_{\text{dBR}} & D_{\text{dBR}} \end{bmatrix} \colon
	      \begin{bmatrix} \DB \\ \cU \end{bmatrix} \to
	 \begin{bmatrix}\DB^{d} \\ \cY \end{bmatrix}
	      $$
	 is an observable coisometric colligation with transfer function
	      equal to $S(z)$:
	      \begin{equation}  \label{dBRreal}
		 S(z) = D_{\text{dBR}} + C_{\text{dBR}}(I_{\DB} - Z(z)
		 A_{\text{dBR}})^{-1} Z(z) B_{\text{dBR}}.
	       \end{equation}
	      Any other observable, coisometric realization of
	      $S$ is unitarily equivalent to this functional-model realization
	      of $S$.
	     \end{theorem}

\begin{proof}
As observed in Proposition \ref{P:H(KS)}, $\cH(K_{S})$ is invariant under
$S_{j}^{*}$ for each $j = 1, \dots, d$.
  From the easily checked intertwining relations
\begin{equation}  \label{tau-intertwine}
     (S^{L}_{j})^{*} \tau =  \tau (S^{R}_{j})^{*} \text{ for } j = 1, \dots,
     d,
\end{equation}
the fact that $\cH(K_{S})$ is invariant under each $(S^{R}_{j})^{*}$
implies that $ \DB$ is invariant under each
$(S^{L}_{j})^{*}$ for $j = 1, \dots, d$. Hence the formula for
$A_{\text{dBR},j}$ in \eqref{dBRops} defines an operator on
$ \DB$.
The first formula for $B_{\text{dBR},j}$ in \eqref{dBRops}
defines an operator from $\cU$
into $\DB$ by part (3) of Proposition \ref{P:H(KS)}; this
is consistent with the second formula as a consequence of
\eqref{tau-intertwine}.
   From \eqref{model-obs} it follows that the pair $(E, {\mathbf
   S}^{*})$ is observable and therefore, since $C$ and ${\mathbf A}$
   are restrictions of $E$ and ${\mathbf S}$ respectively, the pair
   $(C, {\mathbf A})$ is also observable.
Hence, for $u \in \cU$,  making use of \eqref{model-obs}
gives
$$ C_{\text{dBR}} \bA_{\text{dBR}}^{*v} B_{\text{dBR},j} u
= E ({\mathbf S}^{L})^{*v} \tau S_{j}^{*} (M_{S} \cdot u)
= s_{v \cdot  j} u
$$
and it follows that
\begin{align*}
      D_{\text{dBR}} + C_{\text{dBR}} (I - Z(z) \bA_{\text{dBR}})^{-1} Z(z)
	    B_{\text{dBR}} & = s_{\emptyset} +
	    \sum_{j=1}^{\infty} \sum_{v \in \cF_{d}}
	    C_{\text{dBR}} \bA_{\text{dBR}}^{v} B_{\text{dBR},j} z^{v} z_{j}
	    \\  &
	    = s_{\emptyset} + \sum_{j=1}^{d} \sum_{v \in \cF_{d}} s_{v \cdot
	    j} z^{v} z_{j} = S(z)
	   \end{align*}
and \eqref{dBRreal} follows.

By Proposition \ref{P:H(KS)} we know that
$\cH(K_{S})$ is contractively included in $H^{2}_{\cY}(\cF_{d})$,
is invariant under the backward-shift operators $(S^{R}_{j})^{*}$ given by
\eqref{bs} for $j = 1, \dots, d$ with the difference-quotient
inequality \eqref{DQineq-HKS} satisfied.
Hence, by part (4) of Theorem 2.8 in \cite{BBF1}, it follows that
the kernels $K_{S}$ and $K_{C_{\text{dBR}}, \bA_{\text{dBR}}}$ match:
\begin{equation}  \label{=kernel-dBR}
     K_{S}(z,w) = K_{C_{\text{dBR}}, \bA_{\text{dBR}}}(z,w).
\end{equation}
The fact that ${\mathbf U}_{\text{dBR}}$ is coisometric now follows
from Corollary \ref{C:NC-ADR}.
Finally, the uniqueness statement in Theorem \ref{T:NC3} follows from
Corollary \ref{C:unique-coisometric}.
\end{proof}

\begin{remark} \label{R:dBRreal} {\em The proof of Theorem
      \ref{T:NC3} assumed knowledge of the candidate operators
      \eqref{dBRops} for a realization of $S$ and then amounted to a
      check that these operators work.  We remark here that, once
      $A_{\text{dBR}}$ and $B_{\text{dBR}}$  are chosen so that
      \eqref{=kernel-dBR} holds, one can then solve for
      $B_{\text{dBR},1} \dots B_{\text{dBR},d}$ according to the
      prescription \eqref{want1} in the proof of Theorem
      \ref{T:CAStoB}:
      $$ B_{\text{dBR}}^{*} Z(w)^{*}(I - A_{\text{dBR}}^{*}
      Z(w)^{*})^{-1}C^{*} = S(w)^{*} - s_{\emptyset}^{*}
      $$
      to arrive at the formula for $B_{\text{dBR},j}$ ($j = 1,
      \dots, d$) in formula \eqref{dBRops}.}
\end{remark}

\begin{remark} {\em  It is possible to make all the ideas and results
      of this paper symmetric with respect to ``left versus right''.
        Then the multiplication operator
      $M_{S}$ given by \eqref{multiplication} is really the {\em left}
      multiplication operator
      $$
      M^{L}_{S} = \sum_{v \in \cF_{d}} s_{\alpha} ({\mathbf S}^{L})^{v}
      \colon f(z) \mapsto  S(z) \cdot f(z).
      $$
      It is natural to define the corresponding {\em right} multiplication
      operator $M^{R}_{S}$ by
      $$
      M^{R}_{S} = \sum_{v \in \cF_{d}} s_{\alpha} ({\mathbf S}^{R})^{v}.
      $$
      In the scalar case $\cU = \cY = {\mathbb C}$ where $f(z) \cdot
      S(z)$ makes sense, we have
      $$
        M^{R}_{S} \colon f(z) \mapsto f(z) \cdot  (\tau\circ S)(z)
      $$
      while in general we have
      $$
       M^{R}_{S} \colon \sum_{v \in \cF_{d}} f_{v} z^{v}
       \mapsto \sum_{v \in \cF_{d}} \left[ \sum_{\alpha, \beta \in
       \cF_{d} \colon \alpha \beta = v} s_{\beta^{\top}} f_{\alpha}
       \right] z^{v}.
      $$
      The Schur-class ${\mathcal S}_{nc, d}(\cU, \cY)$ is really the
      {\em left} Schur class ${\mathcal S}^{L}_{nc, d}(\cU, \cY)$.  The
      {\em right} Schur class ${\mathcal S}^{R}_{nc,d}(\cU, \cY)$
      consists of all formal power series $S(z) = \sum_{v \in \cF_{d}}
      s_{v} z^{v}$ for which the associated {\em right} multiplication
      operator $M^{R}_{S} = \sum_{v \in \cF_{d}} s_{v} ({\mathbf
      S}^{R})^{v}$ has operator norm at most 1.
      The kernel $K_{S}(z,w)$ given by \eqref{KS} is really the {\em
      left} kernel $K^{L}_{S}(z,w)$ given by
      $$
        K_{S}(z,w) = K^{L}_{S}(z,w) = \{[I_{\cY} - M^{L}_{S}
        (M^{L}_{S})^{*}](k_{\text{Sz}}( \cdot, w) )\}(z).
      $$
      It is then natural to define the corresponding {\em right} kernel
      $$
      K^{R}_{S}(z,w) = \{[I_{\cY} - M^{R}_{S}
        (M^{R}_{S})^{*}](k_{\text{Sz}}( \cdot, w) )\}(z).
      $$
      Given an output pair $(C, \bA)$, the observability operator
      $\cO_{C, \bA}$ given by \eqref{ob-op} is really the {\em left}
      observability operator $\cO^{L}_{C, \bA}$ with range space
      invariant under the {\em right} backward-shift operators
      $(S^{R}_{j})^{*}$; the corresponding {\em right}
      observability operator $\cO^{R}_{C, \bA}$ is given by
      $$
       \cO^{R}_{C, \bA} \colon x \mapsto \sum_{\alpha \in \cF_{d}} (C
       \bA^{v^{\top}}x)z^{\alpha} = C (I - Z({\mathbf S}^{R})A)^{-1} x
      $$
      and has range space invariant under the {\em left} backward shifts
      $(S^{L}_{j})^{*}$.  The system \eqref{sys} is really a {\em left}
      noncommutative multidimensional linear system with
      {\em left} transfer function \eqref{transfunc}
      $$
       T_{\Sigma^{L}}(z) = D + C (I - Z({\mathbf S}^{L}) A)^{-1}
Z({\mathbf S}^{L}) B.
      $$
      For a given colligation $\bU = \sbm{A & B \\ C & D}$, there is an
      associated {\rm right} transfer function
      $$ T_{\Sigma^{R}}(z) = D + C (I - Z({\mathbf S}^{R}) A)^{-1}
      Z({\mathbf S}^{R}) B
      $$
      associated with the {\em right} noncommutative multidimensional
      linear system
      \begin{equation}  \label{sysR}
    \Sigma^{R} \colon \left\{ \begin{array}{ccc}
		x( \alpha \cdot 1) & = & A_{1} x(\alpha) + B_{1} u(\alpha)  \\
		\vdots &   & \vdots  \\
		x( \alpha \cdot d) & = & A_{d} x(\alpha) + B_{d} u(\alpha) \\
		y(\alpha) & = & C x(\alpha) + D u(\alpha)
		\end{array} \right.
	     \end{equation}
initialized with $x(\emptyset) = 0$.
With these definitions in place, it is straightforward to formulate
and prove mirror-reflected versions of Theorem \ref{T:NC1},
Proposition \ref{P:BBF1}, Theorem \ref{T:CAtoS}, Theorem
\ref{T:CAStoB} (as well as Theorems \ref{T:shift=int} and
\ref{T:BLhomint} to come below); we leave the details to the reader.
With all this in hand, it is then possible to identify the state-space
$\DB = \tau \circ \cH(K^{L}_{S})$ appearing in Theorem
\ref{T:NC3} as nothing other than
$\cH(K^{R}_{S})$.  Thus, the functional-model realization for a given
$S$ as an element of the {\em left} Schur class ${\mathcal
S}^{L}_{nc, d}(\cU, \cY)$
uses as state space the functional-model space $\cH(K^{R}_{S})$
based on the {\em right} kernel $K^{R}_{S}$ while the realization of
$S$ as a member of the {\em right} Schur-class ${\mathcal S}^{R}_{nc,
d}(\cU, \cY)$ uses
as the state space the functional-model  $\cH(K^{L}_{S})$ based on
the {\em left} kernel $K^{L}_{S}$.  Presumably it is possible to have
an $S$ in the left Schur-class ${\mathcal S}^{L}_{nc, d}(\cU, \cY)$
but not in the right Schur-class ${\mathcal S}^{R}_{nc, d}(\cU, \cY)$
and vice-versa, although we have not worked out an example.
With this interpretation, the functional-model realization in Theorem
\ref{T:NC3} becomes a more canonical extension of the classical
univariate case.}
\end{remark}

     Let us say that $S \in \cS_{nc, d}(\cU, \cY)$ is {\em inner} if the
     multiplication operator
     $$
     M_{S} \colon H^{2}_{\cU}(\cF_{d}) \to
     H^{2}_{\cY}(\cF_{d})
     $$
     is isometric; such multipliers are the representers for
     shift-invariant subspaces in Popes\-cu's Fock-space analogue of the
     Beurling-Lax theorem \cite{PopescuNF1} (see also \cite{BBF1}). 
It is now an
     easy matter to characterize which functional-model realizations as
in Theorem
     \ref{T:NC3} go with inner multipliers.

     \begin{theorem}  \label{T:NC4}
        The Schur-class multiplier $S \in {\mathcal S}_{nc,d}(\cU, \cY)$
        is inner if and only if $S$ has an observable, coisometric
        realization \eqref{NCrealization} such that $\bA = (A_{1}, \dots,
        A_{d})$ is strongly stable (see
        \eqref{stronglystable}).
      \end{theorem}

      \begin{proof}  By Corollary \ref{C:unique-coisometric}, any
          observable, coisometric realization is unitarily equivalent to
          the functional-model realization given in Proposition
          \ref{P:H(KS)}.  Note that $S$ is inner if and only if $I -
          M_{S}M_{S}^{*}$ is an orthogonal projection.  From the
          characterization of $\cH(K_{S})$ in part (2) of Proposition
	\ref{P:HKS-char}, we see that
          this last condition occurs if and only if $\cH(K_{S})$ is
          contained isometrically in $H^{2}_{\cY}(\cF_{d})$.  By part
          (3) of Proposition \ref{P:BBF1}, this in turn is equivalent to
          strong stability of $\bA$, and Theorem \ref{T:NC4} follows.
       \end{proof}

       \section{Shift-invariant subspaces and Beurling-Lax
       representation theorems}  \label{S:BL}

       Suppose that $(\bZ, X)$ is an isometric input pair, i.e., $\bZ =
       (Z_{1}, \dots, Z_{d})$ where each $Z_{j}
       \colon \cX \to \cX$ and $X \colon \cY \to \cX$. We say that the
       input pair $(\bZ,X)$ is {\em input-stable} if the associated
       controllability operator
       $$ \cC_{\bZ,X} \colon \sum_{v \in \cF_{d}} f_{v} z^{v}
         \mapsto \sum_{v \in \cF_{d}} \bZ^{v^{\top}} X f_{v}
       $$
       maps $H^{2}_{\cY}(\cF_{d})$ into $\cX$.  We say that the pair
       $(\bZ, X)$ is {\em exactly controllable} if in addition
       $\cC_{\bZ,X}$ maps $H^{2}_{\cY}(\cF_{d})$ onto $\cX$.  In this
       case the associated controllability gramian
       $$
       {\mathcal G}_{\bZ,X}: = \cC_{\bZ,X}(\cC_{\bZ,X})^{*}
       $$
       is strictly positive-definite on $\cX$.   and is the unique
solution
       $H = {\mathcal G}_{\bZ,X}$
       of the Stein equation
       \begin{equation}  \label{control-Stein}
       H - Z_{1} H Z_{1}^{*} - \cdots - Z_{d} H Z_{d}^{*} = X X^{*}.
       \end{equation}
       By considering the similar pair
       $$ (\bZ', X') \text{ with } \bZ' = (Z'_{1}, \dots, Z'_{d}) \text{
       where } Z_{j}' = H^{-1/2} Z_{j} H^{1/2} \text{ and } X' =
       H^{-1/2}X,
       $$
       without loss of generality we may assume that the input pair
       $(\bZ, X)$ is {\em isometric}, i.e., \eqref{control-Stein} is
       satisfied with $H = I_{\cX}$.  We are interested in the case when
       in addition $\bZ^{*}$ is {\em strongly stable} in the sense of
       \eqref{stronglystable}; in this case ${\mathcal G}_{\bZ,X}$ is
the unique
       solution of the Stein equation \eqref{control-Stein}.  We remark
       that all these statements are dual to the analogous statements
       made for observability operators $\cO_{C, \bA}$ since the adjoint
       $(C, \bA): = (X^{*}, \bZ^{*})$ of any input pair $(\bZ, X)$ is an
       output pair.

       Given any isometric input pair $(\bZ,X)$ with $\bZ^{*}$ strongly
       stable, we define a {\em left functional calculus with operator
       argument} as follows.  Given $f \in H^{2}_{\cY}(\cF_{d})$ of the
       form $f(z) = \sum_{v \in \cF_{d}} f_{v} z^{v}$, define
       $$
        (X f)^{\wedge L}(\bZ) = \sum_{v \in \cF_{d}} \bZ^{v^{\top}} X
        f_{v} =: \cC_{\bZ,X} f.
      $$
      We define a subspace $\cM_{\bZ,X}$ to be the set of all solutions
      of the associated homogeneous interpolation condition:
      $$
      \cM_{\bZ,X}:= \{ f \in H^{2}_{\cY}(\cF_{d}) \colon (Xf)^{\wedge
      L}(\bZ) = 0\}.
      $$
      That $\cM_{\bZ,X}$ is invariant under the (right) shift operator $S_{j}$
      follows from the intertwining property $\cC_{\bZ,X} S_{j} = Z_{j}
      \cC_{\bZ,X}$ verified by the following computation:
      \begin{align*}
         \cC_{\bZ,X} S_{j} f& = (X S_{j}f)^{\wedge L}(\bZ)
           = \sum_{v \in \cF_{d}} \bZ^{(v j)^{\top}} X f_{v}
           = Z_{j} \cdot \sum_{v \in \cF_{d}}
          \bZ^{v^{\top}} X f_{v} \\
         & = Z_{j} \cdot (X f)^{\wedge L}(\bZ) = Z_{j} \cC_{\bZ,X} f.
       \end{align*}
       It is easily checked that $\cM_{\bZ,X}$ is closed in the metric
       of $H^{2}_{\cY}({\mathcal F}_{d})$.  Hence, by
       Popescu's Beurling-lax theorem for the Fock space (see
       \cite{PopescuNF1}) it is guaranteed that
       $\cM_{\bZ,X}$ has a representation of the form
       $$
       \cM_{\bZ,X} =
       \theta \cdot H^{2}_{\cU}(\cF_{d}) = \operatorname{Ran}\,
M_{\theta}
       $$
       for an inner multiplier $\theta \in {\mathcal S}_{nc, d}(\cU,
       \cY)$.  Our goal is to understand how to compute a
       transfer-function realization for $\theta$ directly from the
       homogeneous interpolation data $(\bZ, X)$.  First, however, we
       show that shift-invariant subspaces $\cM \subset
       H^{2}_{\cY}(\cF_{d})$ of the form $\cM = \cM_{\bZ,X}$ for an
       admissible input pair $(\bZ, X)$ as above are not as special as
       may at first appear.

       \begin{theorem} \label{T:shift=int}
	Suppose that $\cM$ is a closed, shift-invariant
	subspace of $H^{2}_{\cY}(\cF_{d})$.  Then there is an
	isometric input-pair $(\bZ, X)$ with $\bZ^{*}$ strongly
	stable so that $\cM = \cM_{\bZ, X}$.
       \end{theorem}

       \begin{proof}  If $\cM$ is invariant for the operators $S_{j}$,
	then $\cM^{\perp}$ is invariant for the operators $S_{j}^{*}$
	for each $j = 1, \dots, d$. Hence by Theorem 2.8 from
	\cite{BBF1} there is an observable, contractive output pair $(C,
\bA)$ so
	that $\cM^{\perp} = \cH(K_{C,\bA}) = \operatorname{Ran}\,
	\cO_{C, \bA}$ isometrically.  As $\cM^{\perp} \subset
	H^{2}_{\cY}(\cF_{d})$ isometrically, Proposition \ref{P:BBF1}
	tells us that we may take $(C, \bA)$ isometric and that $\bA$
	is strongly stable.  Let $(\bZ, X)$ be the input pair
	$(\bZ,X) = (\bA^{*}, C^{*})$.  As $\cM^{\perp} =
	\operatorname{Ran}\, \cO_{C, \bA}$, we may compute $\cM$ as
	\begin{align*}
	    \cM = \left( \operatorname{Ran}\, \cO_{C, \bA}
\right)^{\perp} =
	    \operatorname{Ker}\, (\cO_{C, \bA})^{*} =
	    \operatorname{Ker}\, \cC_{\bA^{*}, C^{*}} =
	    \operatorname{Ker}\, \cC_{\bZ, X}
	\end{align*}
	and Theorem \ref{T:shift=int} follows.
	\end{proof}

We now suppose that a shift-invariant subspace is given in the form
$\cM = \cM_{\bZ, X}$ for an admissible homogeneous interpolation data
set and we construct a realization for the associated Beurling-Lax
representer.

\begin{theorem}  \label{T:BLhomint}  Suppose that $(\bZ, X)$ is an
       admissible homogeneous interpolation data set and $\cM_{\bZ, X} =
       \operatorname{Ker}\, \cC_{\bZ,X}$ is the associated
       shift-invariant subspace.  Let $(C, \bA)$ be the output pair
       defined by
       $$
         (C, \bA) = (X^{*}, \bZ^{*})
       $$
       and choose an input space $\cU$ with
       $\operatorname{dim}\, \cU = \operatorname{rank}\, (I_{\cX^{d}
       \oplus \cY} - \sbm{ A \\ C } \sbm{A^{*} & C^{*}})$ and define an
       operator $\sbm{B \\ D } \colon \cU \to \cX^{d} \oplus \cY$ as a
       solution of the Cholesky factorization problem
       $$
        \begin{bmatrix} B \\ D \end{bmatrix}
	  \begin{bmatrix} B^{*} & D^{*} \end{bmatrix} =
	 I_{\cX^{d} \oplus \cY} - \begin{bmatrix} A \\ C
        \end{bmatrix} \begin{bmatrix} A^{*} & C^{*} \end{bmatrix}.
       $$
       Set ${\mathbf U} = \sbm{ A & B \\ C & D}$ and let $\theta \in
       {\mathcal S}_{nc,d}(\cU, \cY)$ be the transfer function of
       ${\mathbf U}$:
       $$
         \theta(z) = D + C (I - Z(z) A)^{-1} Z(z) B.
       $$
       Then $\theta$ is inner and $\cM_{\bZ,X} = \theta \cdot
       H^{2}_{\cU}(\cF_{d})$.
       \end{theorem}

        \begin{proof}  If $(\bZ, X)$ is an
       admissible homogeneous interpolation data set, then $(\bZ, X)$
       is controllable and $\bZ^{*}$ is strongly stable. Since  $(C,
\bA) = (X^{*}, \bZ^{*})$,
       we have $(C, \bA)$ is observable and $\bA$ is strongly stable.
>From the construction
       of ${\mathbf U}$, we know ${\mathbf U}$ is coisometric.  Then by
Theorem \ref{T:NC4}, $\theta$ is inner and hence $I- M_{\theta}
M_{\theta}^{*}$ is the orthogonal projection of
$H^{2}_{\cY}(\cF_{d})$ onto $(\operatorname{Ran} \, M_{\theta})^{\perp}$.
     From part (2) of Proposition \eqref{P:HKS-char} it then follows that
\begin{equation}  \label{HKtheta1}
       \cH(K_{\theta}) = H^{2}_{\cY} \ominus \theta \cdot
       H^{2}_{\cU}(\cF_{d}) \text{ isometrically.}
\end{equation}
On the other hand, again since ${\mathbf U}$ is coisometric, from
Corollary \ref{C:NC-ADR} we see that $K_{\theta} = K_{C, \bA}$ and
hence $\cH(K_{\theta}) = \cH(K_{C, \bA})$.  Since $\bA$ is strongly
stable, Proposition \ref{P:BBF1} tells us that $\cH(K_{C, \bA})$ is
isometrically included in $H^{2}_{\cY}(\cF_{d})$ and is characterized
by
\begin{equation}  \label{HKtheta2}
       \cH(K_{\theta}) = \cH(K_{C, \bA}) = \operatorname{Ran}\, \cO_{C,
       \bA} = \operatorname{Ran}\, (\cC_{\bZ,X)})^{*}.
\end{equation}
Comparing \eqref{HKtheta1} with \eqref{HKtheta2} and taking
orthogonal complements finally leaves us with
$$
\theta \cdot H^{2}_{\cU}(\cF_{d}) = (\operatorname{Ran}\,
(\cC_{\bZ,X})^{*})^{\perp} = \operatorname{Ker}\, \cC_{\bZ,X} =
\cM_{\bZ,X}
$$
and Theorem \ref{T:BLhomint} follows.
\end{proof}


\begin{thebibliography}{99}


\bibitem{AM} J.~Agler and J.~E.~McCarthy,  Complete Nevanlinna-Pick
kernels, {\em J. Funct. Anal.}, {\bf  175} (2000), 111--124.

\bibitem{ADR}
D.~Alpay, A.~Dijksma and J.~Rovnyak, A theorem
         of Beurling-Lax type for Hilbert spaces of functions analytic in
         the unit ball, {\em Integral Equations and Operator Theory}
         \textbf{47} (2003), no.3, 251--274.

\bibitem{A-KV} D.~Alpay and D.S.~Kalyuzhny\u \i-Verbovetzki\u \i,
Matrix-$J$-unitary non-commutative rational formal power series, in
{\em The State Space Method: Generalizations and Applications}
(Ed. D. Alpay and I. Gohberg), pp. 49--113, \textbf{OT 161},
Birkh\"auser, Basel-Boston-Berlin, 2006.


\bibitem{ariaspopescu}
A.~Arias and G.~Popescu, {\em Noncommutative interpolation and Poisson
transforms}, Israel J. Math. {\bf 115} (2000), 205--234.

\bibitem{aron}
N.~Aronszajn, \emph{ Theory of reproducing kernels},
Trans. {A}mer. {M}ath. {S}oc.,  \textbf{68} (1950), 337--404.

\bibitem{arv}
W.~Arveson, Subalgebras of $C^*$-algebras. III. Multivariable
operator theory, {\em Acta Math.} \textbf{181} (1998), no. 2,
159--228.

\bibitem{BB-noncomint} J.A.~Ball and V.~Bolotnikov,
\emph{Interpolation in the noncommutative Schur-Agler class},
preprint.

\bibitem{BBF1}  J.A.~Ball, V.~Bolotnikov and Q.~Fang, Multivariable
backward-shift invariant subspaces and observability operators,
preprint.

\bibitem{BBF2} J.A.~Ball, V.~Bolotnikov and Q.~Fang, Realization and
homogeneous interpolation for multipliers of the Arveson space,
preprint.

\bibitem{BGR} J.A.~Ball, I.~Gohberg and L.~Rodman, {\em Interpolation
of Rational Matrix Functions}, \textbf{OT 45}, Birkh\"auser-Verlag,
Basel-Boston-Berlin, 1990.

\bibitem{BGM1} J.A.~Ball, G.~Groenewald and T.~Malakorn,
Structured noncommutative multidimensional linear systems,
{\em SIAM J. Control \& Optimization} \textbf{44} (2005), 1474--1528.

\bibitem{BGM2} J.A.~Ball, G.~Groenewald and T.~Malakorn, Conservative
structured noncommutative multidimensional linear systems, in
{\em The State Space Method: Generalizations and Applications}
(Ed. D. Alpay and I. Gohberg), pp. 179--223, \textbf{OT 161},
Birkh\"auser, Basel, 2006.

\bibitem{BallRaney} J.A.~Ball and M.~Raney, Discrete-time
dichotomous well-posed linear systems and
generalized Schur-Nevanlinna-Pick interpolation, in preparation.

\bibitem{BTV}
J.A.~Ball, T.~T.~Trent and V.~Vinnikov, Interpolation and
commutant lifting for multipliers on reproducing kernel Hilbert
spaces, in {\em Operator Theory and Analysis:  The M.A.~Kaashoek
Anniversary Volume: Workshop in Amsterdam, Nov.~1997} (Ed. H.~Bart,
I.~Gohberg and A.C.M.~Ran), pp. 89--138, \textbf{OT 122}
Birkh\"auser-Verlag, Basel-Boston-Berlin, 2001.

\bibitem{NFRKHS} J.A.~Ball and V.~Vinnikov, Formal reproducing kernel
Hilbert spaces: the commutative and noncommutative settings, in {\em
Reproducing Kernel Spaces and Applications} (Ed.~D.~Alpay), pp.
77--134, \textbf{OT 143}, Birkh\"auser-Verlag, Basel-Boston, 2003.

\bibitem{Cuntz-scat} J.A.~Ball and V.~Vinnikov, \emph{Lax-Phillips
scattering and conservative linear systems: a Cuntz-algebra
multidimensional setting},  Mem. Amer. Math. Soc.  {\bf 178} (2005), no. 
837.

\bibitem{BES} T.~Bhattacharyya, J.~Eschmeier and J.~Sarkar,
       Characteristic function of a pure commuting contractive tuple,
       {\em Integral Equations and Operator Theory} {\bf 53} (2005), 
       no.1, 23--32.


\bibitem{dbr1}
{L. de} Branges and J.~Rovnyak,
\emph{Canonical models in quantum scattering theory},
in: {\emph Perturbation Theory and its
Applications in Quantum Mechanics} (C.~Wilcox, ed.) pp. 295--392,
{Holt, Rinehart and Winston, New York}, 1966.

\bibitem{dbr2}
L.~de Branges and J.~Rovnyak,
{\emph Square summable power series},
Holt, Rinehart and Winston, New York, 1966.

\bibitem{CJ}  T.~Constantinescu and J.L.~Johnson, \emph{A note on
noncommutative interpolation}, Canadian Math.~Bull. \textbf{46} (2003)
no. 1 59-70.

\bibitem{davidsonpitts}
K.~R.~Davidson and D.~R.~Pitts, \emph{Nevanlinna--Pick interpolation
for non-commutative analytic Toeplitz algebras}, Integral
Equations Operator Theory \textbf{31} (1998), no. 3, 321--337.


\bibitem{HMcCV} J.W.~Helton, S.A.~McCullough and V.~Vinnikov, {\em
Noncommutative convexity arises from linear matrix inequalities}, 
Preprint.

\bibitem{KFA} R.F.~Kalman, P.L.~Falb and M.A.~Arbib, {\em Topics in
Mathematical System Theory}, McGraw-Hill, 1969.

\bibitem{KV-V} D.S.~Kalyuzhny\u \i-Verbovetzki\u{\i} and V.~Vinnikov,
{\em Non-commutative positive kernels and their matrix evaluations},
Proc.~Amer.~Math.~Soc. \textbf{134} (2006) no.~3, 805-816.


\bibitem{mt}
S.~McCullough and T.~T.~Trent, Invariant subspaces and
Nevanlinna-Pick kernel, {\em J.~Funct.~Anal.} {\bf 178} (2000), no. 1,
226--249.

\bibitem{MS-Annalen} P.S.~Muhly and B.~Solel, {\em Hardy algebras,
$W^{*}$-correspondences and interpolation theory}, Math.~Annalen
\textbf{330} (2004), 353-415.

\bibitem{MS-Schur} P.S.~Muhly and B.`Solel, {\em Schr class operator
functions and automorphisms of Hardy algebras}, Preprint.

\bibitem{NF} B.~Sz.-Nagy and C.~Foia\c{s}, {\em Harmonic Analysis of
Operators on Hilbert Space}, North-Holland, Amsterdam-London, 1970.


\bibitem{PopescuNF1}
G.~Popescu, Multi-analytic operators and some factorization
theorems, {\em Indiana Univ.~Math.~J.} \textbf{38} (1989), no.~3,
693-710.

\bibitem{PopescuNF2} G.~Popescu, Characteristic functions for
infinite sequences of noncommuting operators, {\em J.~Operator
Theory} \textbf{22} (1989) no. 1, 51-71.

\bibitem{popescu1}
G.~Popescu, Interpolation problems in several variables,
{\em J.~Math.~Anal.~Appl.}, \textbf{227} (1998), 227--250.

\bibitem{Popescu-Nehari} G.~Popescu, \emph{Multivariable Nehari
problem and interpolation}, J.~Funct. Anal. \textbf{200}
(2003), no.~2, 536--581.

\bibitem{Popescu-Memoir} G.~Popescu, \emph{Entropy and multivariable
interpolation},  Mem. Amer. Math. Soc. {\bf 184} (2006), No.868.

\bibitem{sarasonbook}
D.~Sarason, {\emph Sub-{H}ardy {H}ilbert spaces in the unit disk},
John Wiley and Sons Inc., New York, 1994.



\end{thebibliography}
\end{document}